\documentclass{article}

\author{Dixy Msapato}
\usepackage{amsmath,enumerate,amsthm,amsfonts,amssymb,latexsym,bbm}
\usepackage[super]{nth}
\usepackage[a4paper]{geometry}
\usepackage[hidelinks]{hyperref}
\usepackage[UKenglish]{babel}
\usepackage{float}
\usepackage{dirtytalk}
\usepackage{tikz}
\usepackage{tikz-cd}
\usepackage{tikz-qtree}
\usepackage{mathtools}

\newcommand{\myfrac}
    [2]{\begin{array}{@{}c@{}}#1 \\[-0.75ex]#2\end{array}}
\DeclareRobustCommand{\stirling}{\genfrac\{\}{0pt}{}}
\providecommand{\Ext}{\mathop{\rm Ext}\nolimits}
\newcommand{\module}{\mathop{\rm mod}\nolimits}

\newdimen\R
\newdimen\T
\newcommand\blfootnote[1]{%
  \begingroup
  \renewcommand\thefootnote{}\footnote{#1}%
  \addtocounter{footnote}{-1}%
  \endgroup
}

\theoremstyle{definition}
\newtheorem{theorem}{Theorem}[section]
\newtheorem{definition}[theorem]{Definition}
\newtheorem{proposition}[theorem]{Proposition}

\newtheorem{lemma}[theorem]{Lemma}

\newtheorem{example}[theorem]{Example}

\def\t{\tau}

\title{COUNTING THE NUMBER OF $\t$-EXCEPTIONAL SEQUENCES OVER NAKAYAMA ALGEBRAS}
\date{}

\begin{document}
\maketitle
\begin{abstract}
The notion of a $\t$-exceptional sequence was introduced by Buan and Marsh in \cite{Buan2018aa} as a generalisation of an exceptional sequence for finite dimensional algebras. We calculate the number of complete $\t$-exceptional sequences over certain classes of Nakayama algebras. In some cases, we obtain closed formulas which also count other well known combinatorial objects and exceptional sequences of path algebras of Dynkin quivers. 
\end{abstract}
\tableofcontents
\blfootnote{Keywords: $\t$-Exceptional sequence, Exceptional sequence, Nakayama Algebras, $\t$-perpendicular category, $\t$-rigid module, Restricted Fubini numbers. \\
	Contact: mmdmm@leeds.ac.uk}
\section{Introduction}
Let $A$ be a finite dimensional algebra over a field $\mathbb{F}$, where $\mathbb{F}$ is algebraically closed. Let $\module A$ be the category of finitely generated left $A$-modules. A left $A$-module $M$ is called \textit{exceptional} if Hom$(M,M) \cong \mathbb{F}$ and Ext$^{1}_{A}(M,M)=0$. A sequence of indecomposable modules $(M_1,M_2, \dots, M_r)$ is called an \textit{exceptional sequence} if for each pair $(M_i,M_j)$ with $1 \leq i < j \leq r$, we have that Hom$(M_j,M_i) = \Ext^{1}_{A}(M_j,M_i)=0$ and each $M_k$ is exceptional for $1 \leq k \leq r$. Exceptional sequences were first introduced in the context of algebraic geometry by \cite{bondal1989representation}, \cite{gorodentsev1989exceptional} and \cite{gorodentsev1987exceptional}. 

Exceptional sequences exhibit some interesting behaviours.  It was shown by Crawley-Boevey \cite{crawley1993exceptional} and Ringel \cite{ringel1994braid} that there is a transitive braid group action on the set of exceptional sequences. Igusa and Schiffler give \cite{igusa2010exceptional} a characterisation of exceptional sequences for hereditary algebras using the fact that the product of the corresponding reflections is the inverse Coxeter element of the Weyl group. The exceptional sequences for $\module \mathbb{A}_{r}$, where $\mathbb{A}_{r}$ is the path algebra of a Dynkin type A quiver are classified in \cite{garver2015combinatorics} using combinatorial objects called strand diagrams. The exceptional sequences over path algebras of type A, were also characterised using non-crossing spanning trees in \cite{araya2013exceptional}.  A natural question for exceptional sequences is to ask how many there are. The number of them has been computed for all the Dynkin algebras in \cite{seidel2001exceptional} and \cite{obaid2013number}. 

Exceptional sequences have been subject to a number of generalisations. Igusa and Todorov introduce the \textit{signed exceptional sequences} in \cite{igusa2017signed}. More recently, \textit{weak exceptional sequences} were introduced and studied by Sen in \cite{sen2019weak}. Finally, Buan and Marsh introduce in \cite{Buan2018aa} the \textit{signed $\t$-exceptional sequences} and \textit{$\t$-exceptional sequences}. It is $\t$-exceptional sequences which are the subject of this paper. An $A$-module $M$ is called \textit{$\t$-rigid} if Hom$(M,\t M)=0$, see Definition 0.1  in \cite{AdachiIyamaReiten}.  The \textit{$\t$-perpendicular category} of $M$ in $\module A$ is the subcategory $J(M)=  M^{\perp} \cap {^{\perp}( \t M)}$, see Definition 3.3 in \cite{Jasso}.  A sequence of indecomposable modules $(M_1,M_2, \dots, M_r)$ in $\module A$ is called a $\t$-exceptional sequence if $M_r$ is $\t$-rigid in $\module A$ and $(M_1, M_2, \dots , M_{r-1})$ is a $\t$-exceptional sequence in $J(M_r)$. 

Our main results are derivations of closed formulas for the number of complete $\t$-exceptional sequences in the module categories of certain Nakayama algebras. Most notably, we see that the complete $\t$-exceptional sequences of the linear radical square zero Nakayama Algebras $\Gamma_{n}^{2}$ are counted by the \textit{restricted Fubini numbers} $F_{n,\leqslant 2}$ \cite{Mezo}. The numbers $F_{n,\leqslant 2}$ count the number of ordered set partitions of the set $\{1,2, \dots, n\}$ with blocks of size at most two. In the case for the cyclic Nakayama algebra $\Lambda_{n}^{n}$, we get that the complete $\t$-exceptional sequences are counted by the sequence $n^{n}$. We remark that this sequence also counts the number of complete exceptional sequences for the Dykin algebras of quivers of type B and C, as shown in \cite{obaid2013number}.

\section{Definitions and Notation}


Let $A$ be a basic finite dimensional algebra over a field $\mathbb{F}$ which is algebraically closed. Let $\module A$ be the category of finite dimensional left $A$-modules. Denote by $\mathcal{P}(A)$ the full subcategory of projective objects in $\module A$. If $\mathcal{T}$ is a subcategory of $\module A$, we say an $A$-module $M$ in $\mathcal{T}$ is \textit{$\Ext$-projective} in $\mathcal{T}$ if $\text{Ext}_{A}^{1}(M,\mathcal{T})=0$; that is to say $\text{Ext}_{A}^{1}(M,T)=0$ for all $T \in \mathcal{T}.$ We will then write $\mathcal{P}(\mathcal{T})$ to denote the direct sum of indecomposable the Ext-projective modules in $\mathcal{T}$. In everything that follows, we make the assumption that all subcategories are full and closed under isomorphism. We will also take all objects to be basic where possible and they  will be considered up to isomorphism. 

For an additive category $\mathcal{C}$ and an object $X$ in $\mathcal{C}$, we denote by add $X$ the additive subcategory of $\mathcal{C}$ generated by $X$. This is the subcategory of $\mathcal{C}$ with objects the direct summands of direct sums of copies of $X$. For a subcategory $\mathcal{X} \subseteq \mathcal{C}$, we define $^{\perp}\mathcal{X} := \{ Y \in \mathcal{C} : \text{Hom}(Y,X) =0 \text{ for all } X \in \mathcal{X} \}$ and we similarly define $\mathcal{X}^{\perp}$. If $\mathcal{C}$ is skeletally small and Krull-Schmidt, we denote by ind($\mathcal{C}$) the set of isomorphism classes of indecomposable objects in $\mathcal{C}$. For any basic object $X$ in $\mathcal{C}$, let $\delta(X)$ denote the number of indecomposable direct summands of $X$. We fix $\delta(A)$ to be $n$, where $n \geq 1$ is a positive integer. 

Let $\tau$ denote the Auslander-Reitein translate of $\module A$. 
\begin{definition}{\textbf{$\tau$-rigid and $\t$-tilting}}{\cite[Definition 0.1]{AdachiIyamaReiten}}.
A left $A$-module $M$ is said to be \textit{$\tau$-rigid} if Hom($M,\tau M$)=0. If furthermore $\delta(M) = n$, we say that $M$ is \textit{$\t$-tilting}. 
\end{definition}

\begin{definition}{\textbf{$\tau$-perpendicular category}}{\cite[Definition 3.3]{Jasso} .}
Let $M$ be a basic $\t$-rigid left $A$-module. The \textit{$\tau$-perpendicular category} associated to $M$ is the subcategory of $\module A$ given by $J(M) :=  M^{\perp}  \cap {^{\perp}(\tau M}) $.
\end{definition}

\begin{definition}{\textbf{$\tau$-exceptional sequence}}{\cite[Definition 1.3]{Buan2018aa}.}
Let $k$ be a positive integer. A sequence of indecomposable modules $(M_1, M_2, \dots , M_k)$  in $\module A$ is called a \textit{$\tau$-exceptional sequence} in $\module A$ if $M_k$ is $\tau$-rigid in $\module A$ and $(M_1, M_2, \dots , M_{k-1})$ is a $\tau$-exceptional sequence in $J(M_k)$. If $k=n$ we say that the sequence is a complete $\t$-exceptional sequence.
\end{definition}

Let $Q$ be a finite quiver on $n$ vertices labelled by the set $\{1, 2, \dots, n\}$. A path $p$ in $Q$ from the vertex $v_1$ to the vertex $v_m$ is a sequence of vertices $p = (v_1, v_2, v_3, \dots,v_{m-1}, v_m)$ such that $(v_{j},v_{j+1})$ is an arrow in $Q$ for all $1 \leq j \leq m-1.$ The positive integer $m$ is called the length of $p$ and it is denoted by $l(p)$. The path algebra $\mathbb{F}Q$ of the quiver $Q$ is the $\mathbb{F}$-algebra with basis all paths of $Q$ and multiplication defined by concatenation of paths. The the arrow ideal $R_Q$ of $\mathbb{F}Q$ is defined to be the two-sided ideal generated by all arrows in $Q$. The arrow ideal has vector space decomposition given by, 
$$R_Q = \mathbb{F}Q_1 \oplus \mathbb{F}Q_2 \oplus \dots \oplus \mathbb{F}Q_l \oplus \dots$$
where $\mathbb{F}Q_l$ is the subspace of $\mathbb{F}Q$ with basis the set $Q_l$ of paths of length $l$. The $l^{th}$ power of the arrow ideal is $$R_{Q}^{l} := \bigoplus_{m \geq l}\mathbb{F}Q_{m},$$
and it has a basis consisting of all paths of length greater than or equal to $l$. 

For a positive integer $n \geq 1$, let $A_n$ denote the linearly oriented quiver with $n$ vertices, 

\begin{center}
\begin{tikzcd}
    1 \arrow[r,"\alpha_{1}"]& 2 \arrow[r,"\alpha_{2}"] & 3 \arrow[r, "\alpha_{3}"] & \dots \arrow[r,"\alpha_{n-2}"] & n-1 \arrow[r,"\alpha_{n-1}"] & n.
\end{tikzcd}
\end{center}

Let $C_n$ be the linearly oriented $n$-cycle.

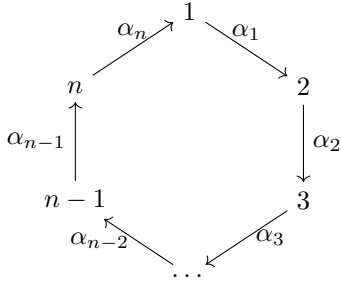
\begin{center}
\begin{tikzpicture}
\node (dots) at (0,0) {$\dots$};
\node (3) at (1.5,1) {$3$};
\node (2) at (1.5,2.5) {$2$};
\node (1) at (0,3.5) {$1$};
\node (n) at (-1.5,2.5) {$n$};
\node (n-1) at (-1.5,1) {$n-1$};

\draw[->] (n-1) -- node[midway, left]{$\alpha_{n-1}$} (n);
\draw[->] (dots) --  node[midway, left]{$\alpha_{n-2}$} (n-1);
\draw[->] (n) --node[midway,above]{$\alpha_{n}$}  (1);
\draw[->] (3) --node[midway,right]{$\alpha_{3}$}  (dots);
\draw[->] (2) --node[midway, right]{$\alpha_{2}$} (3);
\draw[->] (1) --node[midway, above]{$\alpha_{1}$}  (2);
\end{tikzpicture}
\end{center}

  We denote by $\Gamma_{n}^{t}$ the Nakayama algebra $\mathbb{F}A_{n}/R_{Q}^{t}$ and by $\Lambda_{n}^{t}$ the self injective Nakayama algebra $\mathbb{F}C_{n}/R_{Q}^{t}$ where $2 \leq t \leq n$. Throughout the text, we will write $P_{i}$ for the indecomposable projective module at vertex $i$ of the underlying quiver of the algebra $A$ in question. Likewise we will write $S_{i}$ for the simple $A$-module at vertex $i$.

\begin{definition}{\cite[Definition 3.1]{ Schiffler}.}
Let $Q$ be a finite quiver. 
\begin{enumerate}
\item Two paths $p=(v_1, v_2, \dots, v_m)$ and $p^{\prime}=(v_{1}^{\prime}, v_{2}^{\prime}, \dots, v_{m^{\prime}}^{\prime})$ in $Q$ are called \textit{parallel} if $v_1 = v_{1}^{\prime}$ and $v_m = v_{m^{\prime}}^{\prime}$.
\item A \textit{relation} $\rho$ in $Q$ is an $\mathbb{F}$-linear combination $\rho = \sum_{c} \lambda_{c} c $ of parallel paths with $l(c) \geq 2.$
\end{enumerate}
\end{definition}

For a positive integer $n \geq 1$ we will write $(a)_{n}$ to stand for $a$ modulo $n$. We will also write $[i ,j ]_{n}$ for the set $\{(i)_{n},(i+1)_{n}, \dots ,(j-1)_{n},(j)_{n} \}.$

\section{Preliminary Results}
In this section, we will state and prove results which we will be used in later sections to calculate the number of $\tau$-exceptional sequences over the algebras $\Gamma_{n}^{t}$ and $\Lambda_{n}^{t}$. However, our main results are much more general and they apply to other finite dimensional algebras. For this section, we fix an arbitrary finite dimensional $\mathbb{F}$-algebra $A$.

\begin{proposition}{\label{ExtProj}}{\cite[Theorem 2.10]{AdachiIyamaReiten}.}
Let $M$ be a $\tau$-rigid $A$-module. Then the following holds:
\begin{enumerate}
  \item The module $M$ is Ext-projective in $^{\perp}(\t M)$, which is to say that $M$ is in add($\mathcal{P}(^{\perp}(\tau M))).$
  \item The module $T_M := \mathcal{P}(^{\perp}(\t M))$ is a $\t$-tilting $A$-module. 
\end{enumerate}
\end{proposition}
The $A$-module $T_M$ is called the \textit{Bongartz completion} of $M$ in $\module A$. 
\begin{example}{\label{BongartzExample}}
Let $A$ be the algebra $\Gamma_{3}^{2}$ given by the quiver 
\begin{center}
\begin{tikzcd}
    1 \arrow[r,"\alpha_{1}"]& 2 \arrow[r,"\alpha_{2}"] & 3, & & &
\end{tikzcd}
\end{center}
 subject to the relation $\alpha \beta=0$. The Auslander-Reiten quiver of $\module \Gamma_{3}^{2}$ is as follows, 
\begin{center}
\begin{tikzpicture}

\node (3) at (0,0) {$3$};
\node (2) at (2,0) {$2$};
\node (1) at (4,0) {$1$};
\node (23) at (1,1) {$\myfrac{2}{3}$};
\node (12) at (3,1) {$\myfrac{1}{2}$};

\draw[dashed] (3) -- (2);
\draw[dashed] (2)--(1) ;
\draw[->] (3) -- (23);
\draw[->] (23) -- (2);
\draw[->] (2) -- (12);
\draw[->] (12) -- (1);
    		
\end{tikzpicture}.
\end{center}
For the $\Gamma_{3}^{2}$-module $M=1$, $\text{ind}(^{\perp}(\t 1)) =\text{ind} (^{\perp}2) = \{3,\myfrac{1}{2},1\}$. Therefore it is easy to see that $T_{1} = \mathcal{P}(^{\perp}(\t 1)) = 3 \oplus \myfrac{1}{2} \oplus 1$. It is also easy to observe that $T_{1}$ is indeed a $\t$-tilting $\Gamma_{3}^{2}$-module.
\end{example}

\begin{proposition}{\cite[Lemma 2.1]{AdachiIyamaReiten}.}
Let $I$ be an ideal of $A$, and let $M,N$ be $A/I$-modules. Then we have the following: 
\begin{enumerate}
  \item If $\text{Hom}_{A}(M, \t N)=0$ then $\text{Hom}_{A/I}(M, \t_{A/I} N)=0.$
  \item If $I = \langle e \rangle $ for some idempotent $e \in A$, then $\text{Hom}_{A}(M, \t N) =0$
   if and only if $\text{Hom}_{A/I}(M, \t_{A/I} N)=0.$ 
\end{enumerate}
\end{proposition}

The following lemma is well known and it will be important in this paper. 
\begin{lemma}{\label{mainlemma}}
Let $Q$ be a finite simple quiver with vertex set $\{1,2, \dots, n\}$. Let $I$ be the ideal of $\mathbb{F}Q$ generated by relations on $Q$ where each relation is a path in $Q$ and take $A = \mathbb{F} Q/I$. For some $j \in \{1,2, \dots, n\}$ let $Q^{(j)}$ be the quiver obtained from $Q$ by removing the vertex $j$ and any arrows incident to $j$. Let $I^{(j)} \subset I$ be the ideal of $\mathbb{F}Q$ generated by the generating relations of $I$ defined by paths of $Q$ not containing the vertex $j$ and take $B = \mathbb{F}Q^{(j)}/I^{j}$. Then $B \cong A/ \langle e_j \rangle$ as an $\mathbb{F}$-algebra , where $e_j$ is the idempotent at vertex $j$ of $\mathbb{F}Q$. 
\end{lemma}

\begin{theorem}{\label{mainJasso}}{\cite[Theorem 3.8]{Jasso}.}
Let $A$ be a finite dimensional algebra and $M$ a basic $\t$-rigid $A$-module. Let $T_M$ be the Bongartz completion of $M$ in $ \module A$. Let $E_{M} = \text{End}_{A}(T_M)$ and $D_M =E_{M}/ \langle e_M \rangle$, where $e_M$ is the idempotent corresponding to the projective $E_{M}$-module Hom$_{A}(T_M,M)$. Then there is an additive exact equivalence of categories between the category $J(M)$, (the $\t$-perpendicular category of $M$ in $\module A$) and the category $\module D_M$. Moreover, if $M$ is indecomposable we have that $\delta(D_M) = \delta(A)-1$. 

\end{theorem}

We now prove some results which will be crucial in our strategy for calculating the number of $\t$-exceptional sequences in $\module A$. 
\begin{definition}\textbf{Interleaving.}
Let $X=(X_1, X_2, \dots, X_s)$ and $Y=(Y_1, Y_2, \dots ,Y_t)$ be sequences. An interleaved sequence of $X$ and $Y$ is a sequence $Z = (Z_1 ,Z_2, \dots ,Z_{s+t})$ with $Z_i \in \{X_j : 1 \leq i \leq s\} \cup \{Y_j : 1 \leq j \leq t \}$ such that the subsequence of $Z$ containing only elements $X$ or $Y$ is precisely $X$ or $Y$ respectively. 
\end{definition}

\begin{example}
For example let $X = \left(5, \myfrac{4}{5}, 6\right)$ and $Y=\left(2,\myfrac{1}{2}\right)$ be sequences in $\module \Gamma_{6}^{2}$. The sequence $Z=\left(2, 5, \myfrac{4}{5},\myfrac{1}{2}, 6\right)$ is an interleaved sequence of $X$ and $Y$. However $W=\left(\myfrac{4}{5}, 5,2,6,\myfrac{1}{2}\right)$ is not an interleaved sequence of $X$ and $Y$ because the subsequence containing only elements of $X$ is not equal to $X$. 
\end{example}


Let $A$ and $B$ be finite-dimensional $\mathbb{F}$-algebras and let $\module A$ and $\module B$ be the categories of finitely generated left $A$-modules and left $B$-modules respectively. We may consider the category $\module A \oplus \module B$, the direct product category of $\module A$ and $\module B$. The objects of $\module A \oplus \module B$ are pairs $(M,N)$ with $M \in \module A$ and $N \in \module B$. A morphism between  a pair of objects, $(M_1,N_1) \text{ and } (M_2,N_2)$ in $\module A \oplus \module B$ is  a pair of morphisms $(f : M_1 \rightarrow M_2, g: N_1 \rightarrow N_2)$ where $f \in \module A$ and $g \in \module B$. The indecomposable objects of $\module A \oplus \module B$ are pairs $(M,0)$ and $(0,N)$ where $M$ and $N$ are indecomposable in their respective categories. The category $\module A \oplus \module B$ is an abelian category, in fact, there is an exact, additive equivalence to $\module (A\times B).$ The category $\module A \oplus \module B$ also has an Auslander-Reiten translate $\t_{A, B}$ which acts in the obvious way i.e. $\t_{A, B}(M,0)= (\t_A M, 0)$ and $\t_{A, B}(0,N) = (0,\t_B N)$. It is easy to see that the above exact equivalence preserves the Auslander-Reiten translations, since irreducible morphisms, left minimal almost split and right minimal almost split morphisms are preserved under equivalence of categories. Let $M$ be an $A$-module, we identify $M$ with the object $(M,0)$  in $\module A \oplus \module B$. We like wise identify the $B$-module $N$ with the object $(0,N)$ in $\module A \oplus \module B$. It is easy to observe that $(M,0)$ is $\t$-rigid in $\module A \oplus \module B$ if and only if $M$ is $\t$-rigid in $\module A$. The similar statement for $(0,N)$ and $N$ is also true. 

\begin{theorem}{\label{main1}}
Let $A$ and $B$ be finite dimensional $\mathbb{F}$-algebras. Suppose $X=(X_1, X_2, \dots, X_s)$ is a $\t$-exceptional sequence in $\module A$ and $Y=(Y_1,Y_2, \dots, Y_t)$ is a $\t$-exceptional sequence in $\module B$.  Suppose $Z=(Z_1,Z_2, \dots, Z_{s+t})$ is an interleaved sequence of $X$ and $Y$. Then $Z$ is a $\t$-exceptional sequence in $\module A \oplus \module B$.
\begin{proof} 
We prove this by induction on $s+t$. 
For the base case, suppose $s+t=1$. Without loss of generality suppose $t=0$, so $Z=(X_1)$. By assumption, $X_1$ is $\t$-rigid in $\module A$, so it is $\t$-rigid in $\module A \oplus \module B$. This completes the base case.

Suppose the statement is true for $s+t=m$. We consider the $s+t = m+1$ case. Suppose the sequence $Z=(Z_1, Z_2, \dots, Z_{s+t})$ is an interleaved sequence of $X=(X_1, X_2, \dots, X_{s})$ and $Y=(Y_1, Y_2, \dots, Y_{t})$, where $X$ is a $\t$-exceptional sequence in $\module A$ and $Y$ is a $\t$-exceptional sequence in $\module B$. Suppose without loss of generality that $Z_{m+1}$ is in $X$ i.e. $Z_{m+1} = X_{s}$. To show that $Z$ is a $\t$-exceptional sequence in $\module A \oplus \module B$, we need to show that $Z_{m+1}$ is $\t$-rigid in $\module A \oplus \module B$ and that $(Z_1, Z_2, \dots, Z_{m})$ is a $\t$-exceptional sequence in $J_{(A,B)}(Z_{m+1})$, the $\t$-perpendicular category of $Z_{m+1}$ in $\module A \oplus \module B$. By assumption, $Z_{m+1}$ is $\t$-rigid in $\module A$, so it is $\t$-rigid in $\module A \oplus \module B$. Observe that Hom$_{\module A \oplus \module B}(X_s,N) = \text{Hom}_{\module A \oplus \module B}(N, \t X_{s}) = 0$ for all $N \in \module B$, so it follows that 
$$J_{(A,B)}(Z_{m+1})= \{ U \in \module A \oplus \module B : \text{Hom}_{ \module A \oplus \module B}(X_{s},U)=\text{Hom}_{ \module A \oplus \module B}(U, \t_{A} X_{s})=0 \}$$ 
$$=J_{\module A}(X_s) \oplus \module B,$$
where $J_{\module A}(X_s)$ is the $\t$-perpendicular category of $X_s$ in $\module A$. By theorem $\ref{mainJasso}$, $J_{\module A}(X_s)$ is equivalent to a category of modules over some finite dimensional $\mathbb{F}$-algebra. By assumption, $X$ is a $\t$-exceptional sequence in $\module A$, thus $X^{\prime} = (X_1, X_2, \dots, X_{s-1})$ is a $\t$-exceptional sequence in $J_{\module A}(X_s)$. Moreover, $ Z^{\prime} = (Z_1, Z_2, \dots, Z_{m})$ is an interleaved sequence of $X^{\prime}$ and $Y$, so it follows by the inductive hypothesis that $Z^{\prime}$ is a $\t$-exceptional sequence in $J_{\module A}(X_s) \oplus \module C = J_{(A,B)}(Z_{m+1})$, hence $Z$ is a $\t$-exceptional sequence in $\module A \oplus \module B$. This completes the proof. 
\end{proof}
\end{theorem}
We now prove the converse statement. 

\begin{theorem}{\label{main2}}
Let $A$ and $B$ be finite dimensional $\mathbb{F}$-algebras. Suppose $Z=(Z_1, Z_2, \dots, Z_{m})$ is a $\t$-exceptional sequence in $\module A \oplus \module B$. Then $Z$ is an interleaved sequence of some $X=(X_1, X_2, \dots, X_s)$ and $Y=(Y_1, Y_2, \dots, Y_{m-s})$, such that $X$ is a $\t$-exceptional sequence in $\module A$ and $Y$ is a $\t$-exceptional sequence in $\module B$.
\begin{proof}
We prove this by induction on $m$. 

For the base case, suppose $m=1$, so $Z=(Z_1)$ is a $\t$-exceptional sequence in $\module A \oplus \module B$. The module $Z_1$ either lies in $\module A$ or $\module B$. Suppose without loss of generality that $Z_1 \in \module A$.  So we define the sequence $X :=(Z_1)$ and the sequence $Y$ to be the empty sequence. The sequence $Z$ is trivially an interleaved sequence of $X$ and $Y$. As $Z$ is a $\t$-exceptional sequence in $\module A \oplus \module B$, by definition $Z_1$ is $\t$-rigid in $\module A \oplus \module B$, so $Z_1$ is $\t$-rigid in $\module A$. This completes the base case. 

Now suppose the statement is true for $m=k$. We consider the $m=k+1$ case. The sequence $Z=(Z_1,Z_2, \dots, Z_{k+1})$ is a $\t$-exceptional sequence in $\module A \oplus \module B$, so by definition, $Z_{k+1}$ is $\t$-rigid in $\module A \oplus \module B$ and the sequence $Z^{\prime} = (Z_1,Z_2, \dots, Z_k)$ is a $\t$-exceptional sequence in $J_{(A,B)}(Z_{k+1})$, the $\t$-perpendicular category of $Z_{k+1}$ in $\module A \oplus \module B$. Suppose without loss of generality that $Z_{k+1} \in \module A$. We then observe that Hom$_{\module A \oplus \module B}(Z_{k+1},N) = \text{Hom}_{\module A \oplus \module B}(N, \t Z_{k+1}) = 0$ for all $N \in \module B$, so it follows that 
$$J_{(A,B)}(Z_{m+1})= \{ U \in \module A \oplus \module B : \text{Hom}_{ \module A \oplus \module B}(X_{s},U)=\text{Hom}_{ \module A \oplus \module B}(U, \t_{A} X_{s})=0 \}$$ 
$$=J_{\module A}(Z_{k+1}) \oplus \module B,$$
where $J_{\module A}(Z_{k+1})$ is the $\t$-perpendicular category of $Z_{k+1}$ in $\module A$. By theorem $\ref{mainJasso}$, $J_{\module A}(Z_{k+1})$ is equivalent to a category of modules over some finite dimensional $\mathbb{F}$-algebra. So we may apply the inductive hypothesis to $Z^{\prime}$, hence $Z^{\prime}$ is an interleaved sequence of some $X^{\prime} = (X_1, X_2, \dots, X_s)$ and $Y=(Y_1,Y_2, \dots, Y_{k-s})$, where $X^{\prime}$ is a $\t$-exceptional sequence in $J_{\module A}(Z_{k+1})$ and $Y$ is a $\t$-exceptional sequence in $\module B$.  Since $Z_{k+1}$ is $\t$-rigid in $\module A \oplus \module B$, it is also $\t$-rigid $\module A$, hence $X=(X_1, X_2, \dots, X_s, Z_{k+1})$ is a $\t$-exceptional sequence in $\module A$. Clearly $Z$ is an interleaved sequence $X$ and $Y$, so this completes the proof by induction. 
  
\end{proof}
\end{theorem}

In calculating the $\t$-exceptional sequences, it is very convenient to work with the Auslander-Reiten quivers of the module categories. In the case of $\module \Lambda_{n}^{t}$, we will adopt the approach of \cite{sen2019weak} of identifying the Auslander-Reiten quiver of $\module \Lambda_{n}^{t}$ with the following lattices in $\mathbb{Z}^{2}$. Let $f_1 = (2,0)$ and $f_2 = (1,1)$ in $\mathbb{Z}^{2}$. We identify the Auslander-Reiten quiver of $\module \Lambda_{n}^{t}$ with the lattice,
$$\mathcal{AR}(\Lambda_{n}^{t}) = \{ af_{1} + bf_{2} : 0 \leq b \leq t-1, \text{ with } a,b \in \mathbb{Z} \}. $$

The fundamental domain of the Auslander-Reiten quiver of $\module \Lambda_{n}^{t}$ is identified with the lattice, 
$$\mathcal{F}(\Lambda_{n}^{t}) = \{ af_{1} + bf_{2} : 0 \leq b \leq t-1, 0 \leq a \leq n-1 \text{ with } a,b \in \mathbb{Z} \}. $$

\begin{example}
Here we present the Auslander-Reiten quiver of $\module \Lambda_{3}^{2}$ with it's associated lattice which shows the lattice co-ordinates.
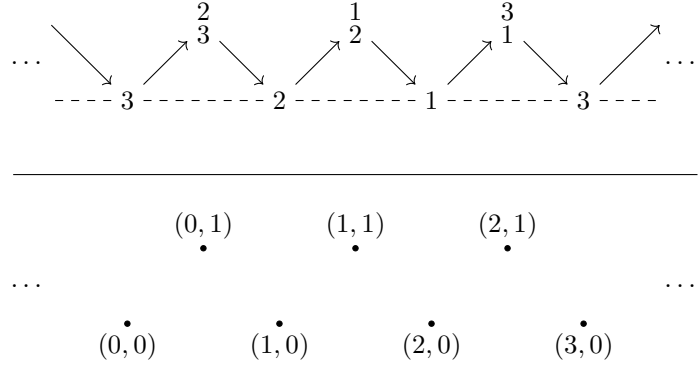
\begin{center}
\begin{tikzpicture}

\node (3) at (0,0) {$3$};
\node (2) at (2,0) {$2$};
\node (1) at (4,0) {$1$};
\node (23) at (1,1) {$\myfrac{2}{3}$};
\node (12) at (3,1) {$\myfrac{1}{2}$};
\node(3s) at (6,0) {$3$};
\node(31) at (5,1) {$\myfrac{3}{1}$};
\node at (7.3,0.5) {$\dots $};
\node at (-1.3,0.5) {$\dots$};

\draw[dashed] (1) -- (3s);
\draw[dashed] (3) -- (2);
\draw[dashed] (2)--(1) ;
\draw[dashed] (3s) -- (7,0);
\draw[dashed] (3) -- (-1,0);
\draw[->] (3) -- (23);
\draw[->] (23) -- (2);
\draw[->] (2) -- (12);
\draw[->] (12) -- (1);
\draw[->] (1)--(31);
\draw[->] (31)--(3s);
\draw[->] (3s) --(7,1);
\draw[->] (-1,1) -- (3);
    		
\end{tikzpicture}
\end{center}

\begin{center}
\begin{tikzpicture}
\draw[solid](1,0) -- (10,0);
\end{tikzpicture}
\end{center}

\begin{center}
\begin{tikzpicture}

\node(3) at (0,-0.3) {$(0,0)$};
\node (2) at (2,-0.3) {$(1,0)$};
\node (1) at (4,-0.3) {$(2,0)$};
\node (3s) at (6,-0.3) {$(3,0)$};
\node at (1,1.3) {$(0,1)$};
\node at (3,1.3) {$(1,1)$};
\node at (5,1.3) {$(2,1)$};
\node at (7.3,0.5) {$\dots $};
\node at (-1.3,0.5) {$\dots$};

\draw[fill=black] (0,0) circle (1pt);
\draw[fill=black] (2,0) circle (1pt);
\draw[fill=black] (4,0) circle (1pt);
\draw[fill=black] (6,0) circle (1pt);
\draw[fill=black] (1,1) circle (1pt);
\draw[fill=black] (3,1) circle (1pt);
\draw[fill=black] (5,1) circle (1pt);
    		
\end{tikzpicture}
\end{center}
\end{example}

We will now recall some standard definitions from \cite{assem2006elements} which we require for the rest of this paper.  Recall that the radical of an $A$-module $M$, denoted by rad$(M)$, is defined to be the intersection of all maximal submodules of $M$. The quotient $M/\text{rad}(M)$ is known as the top of $M$ and is denoted $\text{top}(M)$. The socle of an $A$-module  $M$ denoted soc$(M)$ is the sum of the simple submodules of $M$.

\begin{definition}{\textbf{Radical Series}}{\cite[V.1]{assem2006elements}} Let $M$ be an $A$-module. The radical series of $M$ is defined to be the following sequence of submodules,
$$ 0 \subset \dots \subset \text{rad}^{2}(M) \subset \text{rad}(M) \subset M.$$ 
\end{definition}
Since the left $A$-modules $M$ are finite dimensional as $\mathbb{F}$-vector spaces, there exists a least positive integer $m$ such that rad$^{m}(M)=0.$ The integer $m$ is called the length of the radical series and we denote it by $l(M) = m$. We will also refer to $l(M)$ as the length of the module $M$.

\begin{proposition}{\label{radAssem}}{\cite[V.3.5, V.4.1, V.4.2]{assem2006elements}}
Let $A$ be a basic connected Nakayama algebra and let $M$ be an indecomposable $A$-module. Then there exists some $1 \leq i \leq n$ and $1 \leq j \leq l(P_i)$, such that $M \cong P_{i}/\text{rad}^{j}(P_{i})$ and $j=l(M)$. Moreover, if $M$ is not projective, we have that $\t M \cong \text{rad}(P_{i})/\text{rad}^{j+1}(P_i)$ and $l(\t M) = l(M)$.
\end{proposition}
So we see that modules $M$ of Nakayama algebras are uniquely determined by their top, top$(M)$ and their length $l(M)$.

\begin{proposition}{\label{homAdachi}}{\cite[Lemma 2.4]{Adachi2015tau}}
Let $M = P_{j}/\text{rad}^{l}(P_{j})$ and $N = P_{i}/\text{rad}^{k}(P_{i})$ for $1 \leq i,j,k,l \leq n$. Then the following
conditions are equivalent,
\begin{enumerate}
\item Hom$(M,N) \neq 0$
\item $j \in [i, (i+k-1)]_{n}$ and $(i+k-1)_{n} \in [j , (j+l-1)]_{n}$
\end{enumerate}
\end{proposition}
For every $\Lambda_{n}^{t}$-module $M$, soc($M$) is a simple module. Suppose soc$(M)=S_{i}$, the simple $\Lambda_{n}^{t}$-module associated to the vertex $i$ of $C_{n}$, then the index of the socle of $M$ is defined to be the integer $i$. We write isoc$(M)$ for the index of the socle of $M$.

The  identification of the Auslander-Reiten quiver of $\module \Lambda_{n}^{t}$ with the lattice $\mathcal{AR}(\Lambda_{n}^{t})$ gives rise to the following maps due to \cite{sen2019weak},
$$ L : \module \Lambda_{n}^{t} \rightarrow \mathcal{F}(\Lambda_{n}^{t})$$
$$ M \mapsto (n-\text{isoc}(M)), l(M)-1).$$ 
One may also define $L^{-1}$ from $\mathcal{F}(\Lambda_{n}^{t})$ to $\module \Lambda_{n}^{t}$ in the following way,
$$ L^{-1} : \mathcal{F}(\Lambda_{n}^{t}) \rightarrow \module \Lambda_{n}^{t}$$
$$ (a,b) \mapsto M, \text{ where } l(M)=b+1 \text{ and } \text{isoc}(M)=n-a.$$
In fact the domain of $L^{-1}$ may be extended to $\mathcal{AR}(\Lambda_{n}^{t})$ in the following way; for any $(a,b)$  in $\mathcal{AR}(\Lambda_{n}^{t})$ you may define $L^{-1}(a,b) := L^{-1}((a)_{n},b)$.

\section{The $\Gamma^{2}_{n}$ case}
Let $n \geq 1$ be a positive integer. In this section we will derive a closed formula for the number of complete $\t$-exceptional sequences in $\module \Gamma_{n}^{2}$. Recall that we denote by $A_n$ the linearly oriented quiver with $n$ vertices, 

\begin{center}
\begin{tikzcd}
    1 \arrow[r,"\alpha_{1}"]& 2 \arrow[r,"\alpha_{2}"] & 3 \arrow[r, "\alpha_{3}"] & \dots \arrow[r,"\alpha_{n-2}"] & n-1 \arrow[r,"\alpha_{n-1}"] & n.
\end{tikzcd}
\end{center}

The algebra $\Gamma_{n}^{2}$ is defined to be the $\mathbb{F}$-algebra, $\mathbb{F}A_{n}/R_{Q}^{2}$. This is the path algebra of the quiver $A_{n}$ modulo the relations $\alpha_{i}\alpha_{i+1}=0$ for $1 \leq i \leq n-2$. 

The category $\module \Gamma_{n}^{2}$ has the following Auslander-Reiten quiver.
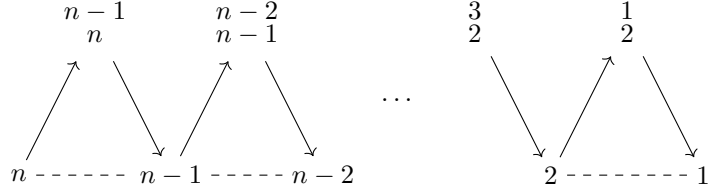
\begin{center}
\begin{tikzpicture}

\node(n) at (0,0) {$n$};
\node (n-1) at (2,0) {$n-1$};
\node (n-2) at (4,0) {$n-2$};
\node (n-1n) at (1,2) {$\myfrac{n-1}{n}$};
\node (n-2n-1) at (3,2) {$\myfrac{n-2}{n-1}$};
\node (32) at (6,2) {$\myfrac{3}{2}$};
\node (2) at (7,0) {$2$};
\node (12) at (8,2) {$\myfrac{1}{2}$};
\node (1) at (9,0) {$1$};
\node (dots) at (5,1) {$\dots$};

\draw[dashed] (n) -- (n-1);
\draw[dashed] (n-1) -- (n-2);
\draw[dashed](2) -- (1);
\draw[->] (n) -- (n-1n);
\draw[->] (n-1n) -- (n-1);
\draw[->] (n-1) -- (n-2n-1);
\draw[->] (n-2n-1) -- (n-2);
\draw[->] (32) -- (2);
\draw[->] (2) -- (12);
\draw[->] (12) -- (1);
    		
\end{tikzpicture}
\end{center}
Our strategy for calculating the number of $\t$-exceptional sequences is straightforward. For each $M$ in ind($\module \Gamma_{n}^{2}$), we will calculate the number of complete $\t$-exceptional sequences ending in $M$. If $M$ is indecomposable, then either $M=P_i$, the projective at vertex $i$ of $A_{n}$, or $M=S_i$, the simple at vertex $i$ of $A_{n}$ (notice that $S_{n}=P_{n}$). In the former case $\t P_{i} = 0$ for $1 \leq i \leq n$ and in the latter case $\t S_j = S_{j+1}$ for $1 \leq j \leq n-1$. In both cases we see that $M$ is $\t$-rigid i.e. every indecomposable $M$ in $\module \Gamma_{n}^{2}$ is $\t$-rigid. We recall that a sequence of indecomposable modules $(M_1, M_2, \dots, M_{n-1},M)$ is a $\t$-exceptional sequence in $\module \Gamma_{n}^{2}$ if $M$ is $\t$-rigid, and $(M_1, M_2, \dots, M_{n-1})$ is a $\t$-exceptional sequence in $J(M)$. Having seen that every indecomposable module $M$ is $\t$-rigid, what is left to do is to calculate $J(M)$ for each indecomposable module. Theorem \ref{mainJasso} and Lemma \ref{mainlemma} are the main tools for these calculations. 

\begin{proposition}{\label{squareproj}}
Let $P_i$ be an indecomposable projective module in $\module \Gamma_{n}^{2}$ for some $1 \leq i \leq n$.  Then the $\t$-perpendicular category of $P_i$ in $\module \Gamma_{n}^{2}$ is $J(P_{i}) \cong \module \Gamma_{i-1}^{2} \oplus \module \Gamma_{n-i}^{2}$.

\begin{proof}
By definition $T_{P_{i}} = \mathcal{P}(^{\perp}(\t P_i))$. Since $\t P_{i} =0$, we have that $^{\perp}(\t P_i)= \module \Gamma_{n}^{2}$. As a result the $\Ext$-projectives of $^{\perp}(\t P_{i})$ are just the projectives of $\module \Gamma_{n}^{2}$, hence $T_{P_i}=\mathcal{P}(^{\perp}(\t P_{i})) = \bigoplus_{j=1}^{n}P_{j}$. Thus the $\mathbb{F}$-algebra $E_{P_{i}} = \text{End}_{\Gamma_{n}^{2}}(T_{P_{i}})$ is precisely given by the path algebra of $A_{n}^{\text{op}},$

\begin{center}
\begin{tikzcd}
    1 & 2 \arrow{l}[swap]{\alpha_{2}} & 3 \arrow{l}[swap]{\alpha_{3}} & \dots \arrow{l}[swap]{\alpha_{4}} & i-1 \arrow{l}[swap]{\alpha_{i-1}} & i \arrow{l}[swap]{\alpha_{i}} & i+1 \arrow{l}[swap]{\alpha_{i+1}} & \dots \arrow{l}[swap]{\alpha_{i+2}} & n-1 \arrow{l}[swap]{\alpha_{n-1}} & n \arrow{l}[swap]{\alpha_{n}}
\end{tikzcd}
\end{center}

 modulo the relations $\alpha_{j}\alpha_{j-1}=0$ for $3 \leq j \leq n$.
Let $A_{n}^{\text{op}(i)}$ be the quiver obtained from $A_{n}^{\text{op}}$ by removing the vertex $i$ and any arrows incident to $i$, 
 
\begin{center}
\begin{tikzcd}
    1 & 2 \arrow{l}[swap]{\alpha_{2}} & 3 \arrow{l}[swap]{\alpha_{3}} & \dots \arrow{l}[swap]{\alpha_{4}} & i-1 \arrow{l}[swap]{\alpha_{i-1}} & & i+1 & \dots \arrow{l}[swap]{\alpha_{i+2}} & n-1 \arrow{l}[swap]{\alpha_{n-1}} & n \arrow{l}[swap]{\alpha_{n}}.
\end{tikzcd}
\end{center} 
 
The quiver $A_{n}^{\text{op}(i)}$ has relations $\alpha_{j}\alpha_{j-1}=0$ for $3 \leq j \leq i-1$ and $i+3 \leq j \leq n$.  By Lemma \ref{mainlemma}, $D_{P_{i}}= E_{P_{i}}/\langle e_{P_{i}} \rangle$ is the path algebra of $A_{n}^{\text{op}(i)}$ modulo its relations. So it follows that $J(P_{i})\cong \module \Gamma_{i-1}^{2} \oplus \module \Gamma_{n-i}^{2}$ by Theorem \ref{mainJasso}.
\end{proof}
\end{proposition}

\begin{proposition}{\label{squaresimp}}
Let $S_{i}$ be a simple non-projective module in $\module \Gamma_{n}^{2}$ for some $1 \leq i \leq n-1.$  Then the $\t$-perpendicular category of $S_{i}$ in $\module \Gamma_{n}^{2}$ is $J(S_{i}) \cong \module \Gamma_{i-1}^{2} \oplus \module \Gamma_{n-i-1}^{2} \oplus \module \Gamma_{1}^{2}$.
\begin{proof}
By definition $T_{S_{i}} = \mathcal{P}(^{\perp}(\t S_{i})).$ Since $S_{i}$ is a simple non-projective indecomposable module $S_{i}$, we have that $\t S_{i} = S_{i+1}.$ Note that the only indecomposable $\Gamma_{n}^{2}$-modules not in $^{\perp}(\t S_{i})$ are $S_{i+1}$ and $P_{i+1}$. Observe also that Ext$_{\Gamma_{n}^{2}}(P_j, {^{\perp}(\t S_{i})})=0$ if $j \neq i+1, 1 \leq j \leq n$. We also have that Ext$_{\Gamma_{n}^{2}}(S_j, {^{\perp}(\t S_{i})}) \neq 0$ for $j \neq i, \text{ and }1 \leq j \leq n$ because $S_{j+1}$ is in $^{\perp}(\t S_{i})$ in these cases. By Proposition \ref{ExtProj}, $S_{i}$ is $\Ext$-projective in $^{\perp}(\t S_{i})$. Therefore $T_{S_{i}} = \mathcal{P}(^{\perp}(\t S_{i})) =  S_{i} \oplus \bigoplus_{j \neq i+1} P_j$.

The $\mathbb{F}$-algebra $E_{S_{i}} = \text{End}_{\Gamma_{n}^{2}}(T_{S_{i}})$ is the path algebra of the following quiver,

\begin{center}
\begin{tikzcd}
    1 & 2 \arrow{l}[swap]{\alpha_{2}} & \dots \arrow{l}[swap]{\alpha_{3}} & i-1 \arrow{l}[swap]{\alpha_{i-1}} & v_{S_{i}} \arrow{l}[swap]{\alpha_{v_{S_{i}}}} & i \arrow{l}[swap]{\alpha_{i}} & & i+2 & \dots \arrow{l}[swap]{\alpha_{i+3}} & n-1 \arrow{l}[swap]{\alpha_{n-1}} & n \arrow{l}[swap]{\alpha_{n}}
\end{tikzcd}
\end{center} 

modulo the relations $\alpha_{j}\alpha_{j-1}=0$ for $ 3 \leq j \leq i-1$ and $i+4 \leq j \leq n$. Here the vertex $v_{S_{i}}$ is the one corresponding to the simple non-projective module $S_{i}$ and the rest correspond to the projective modules $P_{j}$. Consider the following quiver obtained from the one above by removing the vertex $v_{S_{i}}$ and any arrows incident to $v_{S_{i}}$,

\begin{center}
\begin{tikzcd}
    1 & 2 \arrow{l}[swap]{\alpha_{2}} & \dots \arrow{l}[swap]{\alpha_{3}} & i-1 \arrow{l}[swap]{\alpha_{i-1}} & & i & & i+2 & \dots \arrow{l}[swap]{\alpha_{i+3}} & n-1 \arrow{l}[swap]{\alpha_{n-1}} & n \arrow{l}[swap]{\alpha_{n}},
\end{tikzcd}
\end{center} 

it has the relations $\alpha_{j}\alpha_{j-1}=0$ for $ 3 \leq j \leq i-1$ and $i+4 \leq j \leq n$.  By Lemma \ref{mainlemma}, $D_{S_{i}}= E_{S_{i}}/\langle e_{S_{i}} \rangle$ is the path algebra of this quiver modulo its relations. So it follows that $\module D_{S_{i}} \cong \module \Gamma_{i-1}^{2} \oplus \module \Gamma_{n-i-1}^{2} \oplus \module \Gamma_{1}^{2}$. By Theorem \ref{mainJasso}, the statement of this Proposition follows.
\end{proof}
\end{proposition}
Let us denote by $G_n$ the number of complete $\t$-exceptional sequences of $\module \Gamma_{n}^{2}$. When $n=0,1,2$ the $\t$-exceptional sequences coincide with the \say{classical} exceptional sequences since the algebra $\Gamma_{n}^{2}$ is the hereditary Dynkin type $A_n$ algebra $\mathbb{A}_{n}$ in this case. Hence, $G_0=G_1=1$ and $G_2 =3$. 

\begin{lemma}{\label{Gproj}}
Let $P_i$ be the indecomposable projective module in $\module \Gamma_{n}^{2}$ at the vertex $i$ of $A_n$ for some $1 \leq i \leq n$. The number of complete $\t$-exceptional sequences in $\module \Gamma_{n}^{2}$ ending in $P_i$ is, 
$$\binom{n-1}{n-i , i-1}G_{n-i}G_{i-1}.$$
\begin{proof}
Let $(X_1, X_2, \dots,X_{n-1}, P_i)$ be a complete $\t$-exceptional sequence in $\module \Gamma_{n}^{2}$ ending in $P_i$. Then by definition and the fact that $\delta(J(P_i))=n-1$, the sequence $(X_1, X_2, \dots, X_{n-1})$ is a $\t$-exceptional sequence in $J(P_i)$. So to count the number of complete $\t$-exceptional sequences in $\module \Gamma_{n}^{2}$ ending in $P_i$, we just need to count the number of complete $\t$-exceptional sequences in $J(P_i)$. By Lemma \ref{squareproj}, $J(P_i) \cong \module \Gamma_{i-1}^{2} \oplus \module \Gamma_{n-i}^{2}$. By Theorem \ref{main1} and \ref{main2}, the $\t$-exceptional sequences of $J(P_i)$ are interleavings of $\t$-exceptional sequences of $\module \Gamma_{i-1}^{2}$ and $\module \Gamma_{n-i}^{2}$. The number of interleaved sequences coming from a sequence of length $i-1$ and a sequence of length $n-i$ is precisely $\binom{n-1}{n-i,i-1}$. Thus the number of complete $\t$-exceptional sequences ending in $P_i$ is $\binom{n-1}{n-i,i-1}G_{n-i}G_{i-1}$.
\end{proof}
\end{lemma}

\begin{lemma}{\label{Gsimp}}
Let $S_i$ be the indecomposable simple non-projective module in $\module \Gamma_{n}^{2}$ at the vertex $i$ of $A_n$ for some $1 \leq i \leq n-1.$ The number of $\t$-exceptional sequences in $\module \Gamma_{n}^{2}$ ending in $S_i$ is, 
$$\binom{n-1}{n-i-1,i-1}G_{n-i-1}G_{i-1}.$$
\begin{proof}
Let $(X_1, X_2, \dots,X_{n-1}, S_i)$ be a complete $\t$-exceptional sequence in $\module \Gamma_{n}^{2}$ ending in $S_i$. Then by definition and the fact that $\delta(J(S_i))=n-1$, the sequence $(X_1, X_2, \dots, X_{n-1})$ is a complete $\t$-exceptional sequence in $J(S_i)$. Hence to count the number of complete $\t$-exceptional sequences in $\module \Gamma_{n}^{2}$ ending in $S_i$, we just need to count the number of complete $\t$-exceptional sequences in $J(S_i)$. By Lemma \ref{squaresimp}, $J(S_i)\cong \module \Gamma_{i-1}^{2} \oplus \module \Gamma_{n-i-1}^{2} \oplus \module \Gamma_{1}^{2}$. The number of interleaved sequences coming from a sequence of length $i-1$, a sequence of length $n-i-1$ and a sequence of length 1 is precisely $\binom{n-1}{n-i-1,i-1,1}=\binom{n-1}{n-i-1,i-1}$. Thus the number of complete $\t$-exceptional sequences ending in $S_i$ is $\binom{n-1}{n-i-1,i-1}G_{n-i-1}G_{i-1}$.
\end{proof}
\end{lemma}

\begin{theorem}
Let $G_n$ denote the number of complete $\t$-exceptional sequences in $\module \Gamma_{n}^{2}$. Then $G_{n}$ satisfies the recurrence relation, $$G_{n} = \sum_{i=1}^{n}\binom{n-1}{n-i,i-1}G_{n-i}G_{i-1} + \sum_{i=1}^{n-1}\binom{n-1}{n-i-1,i-1}G_{n-i-1}G_{i-1},$$ with initial conditions $G_0=G_1=1.$

\begin{proof}
Let $M$ be an indecomposable in $\module \Gamma_{n}^{2}$, then either $M$ is projective or $M$ simple non-projective. There are $n$ projective indecomposable modules in $\module \Gamma_{n}^{2}$ denoted by $P_i$ for $1 \leq i \leq n.$ There are  $n-1$ simple non-projective indecomposable modules in $\module \Gamma_{n}^{2}$ denoted by $S_i$ for $1 \leq i \leq n-1$. Therefore by Lemma \ref{Gproj} and \ref{Gsimp}, $G_{n} = \sum_{i=1}^{n}\binom{n-1}{n-i,i-1}G_{n-i}G_{i-1} + \sum_{i=1}^{n-1}\binom{n-1}{n-i-1,i-1}G_{n-i-1}G_{i-1}.$
\end{proof}
\end{theorem}
The above Theorem allows us to calculate the first few terms of the sequence $(G_{n})_{n=0}^{\infty}$ as, $$1, 1, 3, 12, 66, 450, 3690, 35280, 385560, 4740120.$$

An ordered set partition of $\{1,2, \dots, n\}$, is a partition of the set $\{1,2, \dots, n \}$ together with a total order on the sets in the partition. We refer to the sets in an ordered partition as blocks. The number $F_{n,\leqslant m}$ counts the number of ordered set partitions of $\{1,2, \dots, n \}$ with blocks of size at most $m$. The sequence $(F_{n,\leqslant m})$ is known as the \textit{restricted Fubini numbers}. For example, when $n=m=2$, the ordered set partitions are $(\{1,2\}), (\{1\},\{2\}) \text{ and } (\{2\},\{1\})$ hence $F_{2,\leqslant 2}=3$. The  \textit{restricted Stirling number of the second kind} denoted by $\stirling{n}{k}_{\leqslant m}$ is the number of (unordered) partitions of $\{1, 2, \dots, n\}$ into $k$ subsets with the restriction that each block contains at most $m$ elements, therefore $$F_{n,\leqslant m} = \sum_{k=0}^{n} k!\stirling{n}{k}_{\leqslant m}.$$

It is shown in \cite[Section 5.4]{Mezo} that the restricted Fubini numbers satisfy the recurrence, 
$$F_{n,\leqslant m} = \sum_{l=1}^{m} \binom{n}{l} F_{n-l,\leqslant m}.$$ 

The sequence $(F_{n,\leqslant 2})$ is listed on the On-line Encyclopedia of Integer Sequences (OEIS) as the sequence A080599. The first terms of this sequence coincide with the first terms we calculated for $(G_n)$ so we would like to prove that it is is the case that $F_{n,\leqslant 2} = G_n$. 

When $m=2$ the recurrence for $F_{n,\leqslant m}$ is given as $F_{n,\leqslant 2} = nF_{n-1,\leqslant 2} + \binom{n}{2}F_{n-2,\leqslant 2}.$ In the paper \cite[Theorem 3.7]{gellert2017degree}, the authors derive the closed formula
$$F_{n,\leqslant 2} = \frac{n!}{\sqrt{3}}\left((\sqrt{3}-1)^{-n-1}-(-\sqrt{3}-1)^{-n-1}\right). $$
An exponential generation function for $F_{n,\leqslant m}$ is given in \cite[Theorem 4]{komatsu2018some}: 
$$\sum_{n=0}^{\infty} F_{n,\leqslant m} \frac{x^{n}}{n!} = \frac{1}{1-x-\frac{x^2}{2!}- \dots \ \frac{x^{m}}{m!}}.$$
We will show that $G_n = F_{n,\leqslant 2}$ by showing that the exponential generation functions  for $G_n$ and $F_{n,\leqslant 2}$ coincide. 

\begin{theorem}{\label{Gmain}}
Let $G_n$ denote the number of complete $\t$-exceptional sequences in $\module \Gamma_{n}^{2}$. The exponential generating function of $G_{n}$ is as follows,
$$\sum_{n=0}^{\infty} G_{n} \frac{x^{n}}{n!} = \frac{1}{1-x-\frac{x^2}{2!}} \text{ . }$$ Therefore $G_n = F_{n,\leqslant 2}$ and $$G_{n} =  \frac{n!}{\sqrt{3}}\left((\sqrt{3}-1)^{-n-1}-(-\sqrt{3}-1)^{-n-1}\right).$$ 
\begin{proof} First let us recall the recurrence relation for $G_n$. 
$$G_{n} = \sum_{i=1}^{n}\binom{n-1}{n-i,i-1}G_{n-i}G_{i-1} + \sum_{i=1}^{n-1}\binom{n-1}{n-i-1,i-1}G_{n-i-1}G_{i-1}.$$
$$= \sum_{i=1}^{n} \frac{n-1}{(n-i)!(i-1)!}G_{n-i}G_{i-1} +  \sum_{i=1}^{n-1} \frac{n-1}{(n-i-1)!(i-1)!}G_{n-i-1}G_{i-1}.$$
Therefore
$$G_{n+1} = \sum_{i=1}^{n+1} \frac{n}{(n+1-i)!(i-1)!}G_{n+1-i}G_{i-1} +  \sum_{i=1}^{n} \frac{n}{(n-i)!(i-1)!}G_{n-i}G_{i-1}.$$
Let $$g(x) = \sum_{n=0}^{\infty} G_{n} \frac{x^{n}}{n!}  \text{ with } g(0)=1,$$ be the exponential generating function of $G_n$. We then have that the first derivative of $g(x)$ is $g^{\prime}(x) = \sum_{n=0}^{\infty} G_{n+1} \frac{x^{n}}{n!}.$ Expanding $G_{n+1}$ in $g^{\prime}(x)$ by the recurrence relation above we obtain the following.
$$g^{\prime}(x) = \sum_{n=0}^{\infty} \left(\sum_{i=1}^{n+1} \frac{n!}{(n+1-i)!(i-1)!}G_{n+1-i}G_{i-1}\right)\frac{x^{n}}{n!} + \sum_{n=0}^{\infty}\left(\sum_{i=1}^{n} \frac{n!}{(n-i)!(i-1)!}G_{n-i}G_{i-1}\right)\frac{x^{n}}{n!}$$
$$= \sum_{n=0}^{\infty} \left(\sum_{i=1}^{n+1} \frac{G_{n+1-i}G_{i-1}}{(n+1-i)!(i-1)!}\right)x^{n} + \sum_{n=0}^{\infty} \left(\sum_{i=1}^{n} \frac{G_{n-i}G_{i-1}}{(n-i)!(i-1)!}\right)x^{n}.$$

Recall the Cauchy product of formal power series is as follows, 
$$\left(\sum_{s=0}^{\infty} a_{s}x^{s}\right)\left(\sum_{t=0}^{\infty} b_{t}x^{t}\right) = \sum_{k=0}^{\infty}c_{k}x^{k} \text{ where } c_{k} = \sum_{l=0}^{k} a_{l}b_{k-l}.$$
By performing a change of variable in $g^{\prime}(x)$ by setting $j = i -1$ and factorising $x$ from the right summand we write,
$$g^{\prime}(x) = \sum_{n=0}^{\infty} \left(\sum_{j=0}^{n}\frac{G_{n-j}G_{j}}{(n-j)!j!}\right)x^{n} +x \sum_{n=0}^{\infty} \left(\sum_{j=0}^{n-1} \frac{G_{n-j-1}G{j}}{(n-j-1)!j!}\right)x^{n-1}.$$
Using the Cauchy product of formal power series, we obtain the following first order non-linear ordinary differential equation. 
$$g^{\prime}(x) = (g(x))^{2} + x(g(x))^{2} = (1+x)(g(x))^2 \text { with initial conditions } g(0)=1.$$
It is easy to check that the unique solution to this ODE is given by, 
$$g(x) = \frac{-2}{-2 + x(x+2)} = \frac{1}{1-x-\frac{x^2}{2}} \text{ . }$$This completes the proof.
\end{proof}
\end{theorem}

\section{The $\Lambda_{n}^{2}$ case}
Let $n \geq 1$ be a positive integer. In this section we will derive a closed formula for the number of complete $\t$-exceptional sequences in $\module \Lambda_{n}^{2}$. Recall that we denote by $C_n$ the linearly oriented $n$-cycle. 
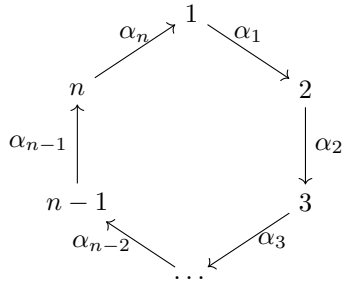
\begin{center}
\begin{tikzpicture}
\node (dots) at (0,0) {$\dots$};
\node (3) at (1.5,1) {$3$};
\node (2) at (1.5,2.5) {$2$};
\node (1) at (0,3.5) {$1$};
\node (n) at (-1.5,2.5) {$n$};
\node (n-1) at (-1.5,1) {$n-1$};

\draw[->] (n-1) -- node[midway, left]{$\alpha_{n-1}$} (n);
\draw[->] (dots) -- node[midway, left]{$\alpha_{n-2}$} (n-1);
\draw[->] (n) --node[midway,above]{$\alpha_{n}$}  (1);
\draw[->] (3) --node[midway,right]{$\alpha_{3}$}  (dots);
\draw[->] (2) --node[midway, right]{$\alpha_{2}$} (3);
\draw[->] (1) --node[midway, above]{$\alpha_{1}$}  (2);
\end{tikzpicture}
\end{center}
The algebra $\Lambda_{n}^{2}$ is defined to be the $\mathbb{F}$-algebra, $\mathbb{F}C_{n}/R_{Q}^{2}$. This is the path algebra of the quiver $C_{n}$ modulo the relations $\alpha_{j}\alpha_{(j+1)_{n}}=0$ for $1 \leq j \leq n$. 

The category $\module \Lambda_{n}^{2}$ has the following Auslander-Reiten quiver.
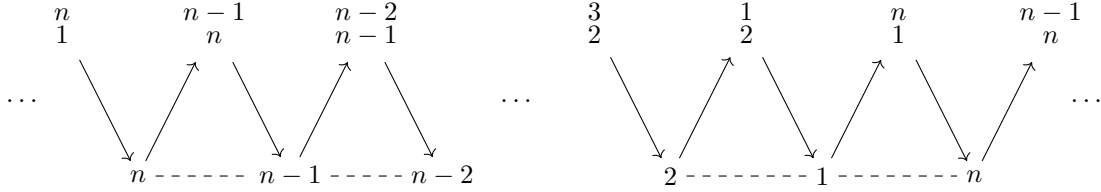
\begin{center}
\begin{tikzpicture}
\node (b) at (-1,2) {$\myfrac{n}{1}$};
\node(fb) at (10,2) {$\myfrac{n}{1}$};
\node (fn) at (11,0) {$n$};
\node(fn-1n) at (12,2) {$\myfrac{n-1}{n}$};
\node(n) at (0,0) {$n$};
\node (n-1) at (2,0) {$n-1$};
\node (n-2) at (4,0) {$n-2$};
\node (n-1n) at (1,2) {$\myfrac{n-1}{n}$};
\node (n-2n-1) at (3,2) {$\myfrac{n-2}{n-1}$};
\node (32) at (6,2) {$\myfrac{3}{2}$};
\node (2) at (7,0) {$2$};
\node (12) at (8,2) {$\myfrac{1}{2}$};
\node (1) at (9,0) {$1$};
\node (dots) at (5,1) {$\dots$};
\node (dota) at (-1.5,1) {$\dots$};
\node (dotas) at (12.5,1) {$\dots$};

\draw[dashed] (1) -- (fn);
\draw[dashed] (n) -- (n-1);
\draw[dashed] (n-1) -- (n-2);
\draw[dashed](2) -- (1);
\draw[->] (n) -- (n-1n);
\draw[->] (n-1n) -- (n-1);
\draw[->] (n-1) -- (n-2n-1);
\draw[->] (n-2n-1) -- (n-2);
\draw[->] (32) -- (2);
\draw[->] (2) -- (12);
\draw[->] (b) -- (n);
\draw[->] (12) -- (1);
\draw[->] (1) -- (fb);
\draw[->] (fb)--(fn);
\draw[->] (fn) -- (fn-1n);
    		
\end{tikzpicture}
\end{center}
We will use the same approach for calculating the number of complete $\t$-exceptional sequences for $\module \Lambda_{n}^{2}$ as we did for $\module \Gamma_{n}^{2}$. If $M$ is indecomposable in $\module \Lambda_{n}^{2}$, then $M=P_i$, the projective at vertex $i$ of $C_{n}$, or $M=S_i$, the simple at vertex $i$ of $C_{n}$. In the former case $\t P_{i} = 0$ and in the latter case $\t S_i = S_{(i+1)_{n}}$. In both cases $M$ is $\t$-rigid i.e. every indecomposable $M$ in $\module \Lambda_{n}^{2}$ is $\t$-rigid. We recall that a sequence of indecomposable modules $(M_1, M_2, \dots, M_{n-1},M)$ is a $\t$-exceptional sequence in $\module \Lambda_{n}^{2}$ if $M$ is $\t$-rigid, and $(M_1, M_2, \dots, M_{n-1})$ is a $\t$-exceptional sequence in $J(M)$. Having seen that every indecomposable module $M$ is $\t$-rigid, what is left to do is to calculate $J(M)$ for each indecomposable module. Theorem \ref{mainJasso} and Lemma \ref{mainlemma} are the main tools for these calculations. 

\begin{proposition}{\label{L2proj}}Let $P_i$ be an indecomposable projective module in $\module \Lambda_{n}^{2}$ for some $1 \leq i \leq n$.  Then the $\t$-perpendicular category of $P_i$ in $\module \Lambda_{n}^{2}$ is $J(P_i) \cong \module \Gamma_{n-1}^{2}$.
\begin{proof}
By definition $T_{P_i} = \mathcal{P}(^{\perp}(\t P_i))$. Since $P_i$ is projective, we have that $\t P_i =0$, therefore $^{\perp}(\t P_i)= \module \Lambda_{n}^{2}$. As a result the $\Ext$-projectives of $^{\perp}(\t P_i)$ are just the projectives of $\module \Lambda_{n}^{2}$, hence $T_{P_i}=\mathcal{P}(^{\perp}(\t P_i)) = \bigoplus_{j=1}^{n}P_{j}$. Thus the $\mathbb{F}$-algebra $E_{P_i} = \text{End}_{\Lambda_{n}^{2}}(T_{P_i})$ is precisely given by the path algebra of the quiver $C_{n}^{\text{op}},$
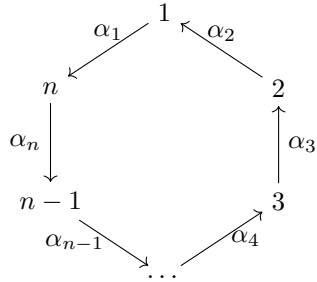
\begin{center}
\begin{tikzpicture}
\node (dots) at (0,0) {$\dots$};
\node (3) at (1.5,1) {$3$};
\node (2) at (1.5,2.5) {$2$};
\node (1) at (0,3.5) {$1$};
\node (n) at (-1.5,2.5) {$n$};
\node (n-1) at (-1.5,1) {$n-1$};

\draw[->] (n) -- node[midway, left]{$\alpha_{n}$} (n-1);
\draw[->] (n-1) --  node[midway, left]{$\alpha_{n-1}$} (dots);
\draw[->] (1) --node[midway,above]{$\alpha_{1}$}  (n);
\draw[->] (dots) --node[midway,right]{$\alpha_{4}$}  (3);
\draw[->] (3) --node[midway, right]{$\alpha_{3}$} (2);
\draw[->] (2) --node[midway, above]{$\alpha_{2}$}  (1);
\end{tikzpicture}
\end{center}
modulo the relations $\alpha_{i}\alpha_{(i-1)_{n}}=0$ for $1 \leq i \leq n$.

Let $C_{n}^{\text{op}(i)}$ be the quiver obtained from $C_{n}^{\text{op}}$ by removing the vertex at $i$ and any arrows incident to $i$.

\begin{center}
\begin{tikzcd}
   i+ 1 & i+2 \arrow{l}[swap]{\alpha_{i+2}} & \dots \arrow{l}[swap]{\alpha_{i+3}} & n-1 \arrow{l}[swap]{\alpha_{n-1}} & n \arrow{l}[swap]{\alpha_{n}} & 1 \arrow{l}[swap]{\alpha_{1}} & \dots \arrow{l}[swap]{\alpha_{2}} & i-2 \arrow{l}[swap]{\alpha_{i-2}} & i-1 \arrow{l}[swap]{\alpha_{i-1}}
\end{tikzcd}
\end{center}
It has the relations $\alpha_{1}\alpha_{n}=0$ and $\alpha_{j}\alpha_{j-1}=0$ for $i+2 \leq j \leq n$ and $2 \leq j \leq i-2$. By Lemma \ref{mainlemma}, $D_{P_i} = E_{P_i}/ \langle e_{P_i} \rangle$ is the path algebra of the quiver $C_{n}^{\text{op}(i)}$ modulo relations.  It is easy to see that in fact $D_{P_{i}}$ is isomorphic to $\Gamma_{n-1}^{2}$. Hence by Theorem \ref{mainJasso}, the $\t$-perpendicular category $J(M) \cong \module \Gamma_{n-1}^{2}$.
\end{proof}
\end{proposition}

\begin{proposition}{\label{L2sim}}
Let ${S_{i}}$ be a simple module in $\module \Lambda_{n}^{2}$ for some $1 \leq i \leq n.$  Then the $\t$-perpendicular category of ${S_{i}}$ in $\module \Lambda_{n}^{2}$ is $J({S_{i}}) \cong \module \Gamma_{n-2}^{2} \oplus \module \Gamma_{1}^{2}.$
\begin{proof}
By definition $T_{S_{i}} = \mathcal{P}(^{\perp}(\t {S_{i}}))$. Since $S_i$ is a simple $\Lambda_{n}^{2}$-module, we have that $\t S_{i} = S_{(i+1)_{n}}$. Note that the only $\Lambda_{n}^{2}$-modules not in $^{\perp} ( \t {S_{i}})$ are $S_{(i+1)_{n}} \text{ and } P_{(i+1)_{n}}.$ Observe also that $\Ext_{\Lambda_{n}^{2}}(P_{j}, {^{\perp}(\t {S_{i}})}) = 0$ for $ j \neq (i+1)_{n}$ and $1 \leq j \leq n$. However for $j \neq i, i+1$, Ext$_{\Lambda_{n}^{2}}(S_{j}, {^{\perp}(\t {S_{i}})}) \neq 0$ because $S_{(j+1)_{n}}$ is in $^{\perp}(\t M)$. By Proposition \ref{ExtProj}, ${S_{i}}$ is $\Ext$-projective in ${^{\perp}(\t {S_{i}})}$.  Hence $T_{S_{i}} = \mathcal{P}(^{\perp}(\t {S_{i}})) = {S_{i}} \oplus \bigoplus_{j \neq (i+1)_{n}} P_{j}$ is the Bongartz completion of ${S_{i}}$.

The $\mathbb{F}$-algebra of $E_{S_{i}} = \text{End}_{\Lambda_{n}^{2}}(T_{S_{i}})$ is given by the path algebra of the quiver, 

\begin{center}
\begin{tikzcd}
    i+2 & i+3 \arrow{l}[swap]{\alpha_{i+3}} & \dots \arrow{l}[swap]{\alpha_{i+4}} & n-1 \arrow{l}[swap]{\alpha_{n-1}} & n \arrow{l}[swap]{\alpha_{n}} & 1 \arrow{l}[swap]{\alpha_{1}} & \dots \arrow{l}[swap]{\alpha_{2}} & i-1 \arrow{l}[swap]{\alpha_{i-1}} & v_{S_{i}} \arrow{l}[swap]{\alpha_{v_{S_{i}}}} & i \arrow{l}[swap]{\alpha{i}}
\end{tikzcd}
\end{center}

modulo the relations $\alpha_{v_{S_{i}}}\alpha_{i-1}=0=\alpha_{1}\alpha_{n} \text{ and } \alpha_{j}\alpha_{j-1}=0$ for $ i+4 \leq j \leq n$ and $2 \leq j \leq i-1$. Here the vertex $v_{{S_{i}}}$ is the one corresponding to the simple module ${S_{i}}$ and the rest correspond to the projective modules $P_{j}$. By Lemma \ref{mainlemma}, $D_{S_{i}} = E_{S_{i}}/ \langle e_{S_{i}} \rangle$ is the path algebra of the quiver obtained from the one above by removing the vertex $v_{S_{i}}$,
\begin{center}
\begin{tikzcd}
    i+2 & i+3 \arrow{l}[swap]{\alpha_{i+3}} & \dots \arrow{l}[swap]{\alpha_{i+4}} & n-1 \arrow{l}[swap]{\alpha_{n-1}} & n \arrow{l}[swap]{\alpha_{n}} & 1 \arrow{l}[swap]{\alpha_{1}} & \dots \arrow{l}[swap]{\alpha_{2}} & i-1 \arrow{l}[swap]{\alpha_{i-1}} &  & i 
\end{tikzcd}
\end{center}
modulo the relations $\alpha_{1}\alpha_{n}=0 \text{ and } \alpha_{j}\alpha_{j-1}=0$ for $ i+4 \leq j \leq n$ and $2 \leq j \leq i-1$. 
So it follows that $\module D_{S_{i}} \cong \module \Gamma_{n-2}^{2} \oplus \module \Gamma_{1}^{2}$. By Theorem \ref{mainJasso} the statement of this Proposition follows.  
\end{proof}
\end{proposition}
Denote by $L_{n}$ the number of complete $\t$-exceptional sequences in $\module \Lambda_{n}^{2}$. 

\begin{theorem}
Let $L_{n}$ be the number of complete $\t$-exceptional sequences in $\module \Lambda_{n}^{2}$. Then $L_{n}$ satisfies the relation, 
$$L_{n} = nG_{n-1} + n(n-1)G_{n-2},$$
with initial conditions $L_{1} =1 \text{ and } L_{2}=4$,  and where $G_{m}$ denotes the number of complete $\t$-exceptional sequences in $\module \Gamma_{m}^{2}$.
\begin{proof}
Suppose $M$ is an indecomposable projective $\Lambda_{n}^{2}$-module, then by Lemma \ref{L2proj}, the $\t$-perpendicular category $J(M) \cong \module \Gamma_{n-1}^{2}$. Suppose $(X_1, X_2, \dots , X_{n-1}, M)$ is a complete $\t$-exceptional sequence ending in $M$ in $\module \Lambda_{n}^{2}$ . Then by the fact that $\delta(J(M))=n-1$ and by definition, the sequence $(X_1, X_2, \dots, X_{n-1})$ is a complete $\t$-exceptional sequence in $J(M) \cong \module \Gamma_{n-1}^{2}$. Hence the number of complete $\t$-exceptional sequences ending in $M$ is $G_{n-1}$, which is the number of complete $\t$-exceptional sequences in $\module \Gamma_{n-1}^{2}$. 

Now suppose $M$ is a simple $\Lambda_{n}^{2}$-module. By Lemma \ref{L2sim},  the $\t$-perpendicular category $J(M) \cong \module \Gamma_{n-2}^{2} \oplus \module \Gamma_{1}^{2}$. Arguing as above the number of complete $\t$-exceptional sequences ending in $M$ is equal to the number of complete $\t$-exceptional sequences in $J(M)$. Since $J(M) \cong \module \Gamma_{n-2}^{2} \oplus \module \Gamma_{1}^{2}$, by Theorem \ref{main1} and \ref{main2}, the $\t$-exceptional sequences of $J(M)$ are interleavings of $\t$-exceptional sequences of $\text{mod}\Gamma_{n-2}^{2}$ and $\text{mod}\Gamma_{1}^{2}$. The number of interleaved sequences coming from a sequence of length $n-2$ and a sequence of length $1$ is precisely $\binom{n-1}{n-2,1}=(n-1)$. Thus the number of complete $\t$-exceptional sequences ending in $M$ is $(n-1)G_{n-2}G_{1} = (n-1)G_{n-2}$. 

An arbitrary indecomposable $\Lambda_{n}^{2}$-module is either projective or simple. There are $n$ projective modules and $n$ simple modules up to isomorphism in $\module \Lambda_{n}^{2}$, hence the number of complete $\t$-exceptional sequences in $\module \Lambda_{n}^{2}$ is $L_{n} =  nG_{n-1} + n(n-1)G_{n-2}$. It then follows easily that $L_{1}=1$ and $L_{2}=4$. 
\end{proof}
\end{theorem}

In the previous section we found the exponential generating function and closed formula for $G_{n}$. Using the above theorem, we can immediately do the same for $L_{n}$.

\begin{theorem}
Let $L_{n}$ denote the number of complete $\t$-exceptional sequences in $\module \Lambda_{n}^{2}$. The exponential generating function of $L_{n}$ is as follows,
$$\sum_{n=0}^{\infty} L_{n} \frac{x^{n}}{n!} = \frac{x+x^{2}}{1-x-\frac{x^2}{2}} \text{ . }$$ 
\begin{proof} Let $h(x) = \sum_{n=0}^{\infty} L_{n} \frac{x^{n}}{n!}$ be the exponential generating function of $L_{n}.$ Let $g(x) = \sum_{n=0}^{\infty} G_{n} \frac{x^{n}}{n!}$ be the exponential generating function of $G_{n}$.
We then recall the recurrence relation of $L_{n}$, 
$$ L_{n} =  nG_{n-1} + n(n-1)G_{n-2}.$$
Therefore the exponential generating function of $L_{n}$ is, 
$$\sum_{n=0}^{\infty} L_{n} \frac{x^{n}}{n!} = \sum_{n=0}^{\infty} nG_{n-1}\frac{x^{n}}{n!} + \sum_{n=0}^{\infty} n(n-1)G_{n-2} \frac{x^{n}}{n!}.$$

$$= \sum_{n=0}^{\infty} G_{n-1} \frac{x^{n}}{(n-1)!} + \sum_{n=0}^{\infty} G_{n-2} \frac{x^{n}}{(n-2)!} $$
$$= x\sum_{n=0}^{\infty} G_{n-1} \frac{x^{n-1}}{(n-1)!} + x^{2}\sum_{n=0}^{\infty} G_{n-2} \frac{x^{n-2}}{(n-2)!}.$$
Therefore 
$$h(x) = xg(x) + x^{2}(g(x)) = (x+x^{2})g(x).$$
By Theorem \ref{Gmain}, $$g(x) = \frac{1}{1-x-\frac{x^{2}}{2}},$$ hence $$h(x) = \frac{x+x^{2}}{1-x-\frac{x^2}{2!}} \text{ . }$$
\end{proof}
\end{theorem}

\begin{theorem}Let $L_{n}$ denote the number of complete $\t$-exceptional sequences in $\module \Lambda_{n}^{2}$. Then $L_{n}$ is given by the closed formula, $$L_{n} = \frac{n!}{\sqrt{3}}\left((\sqrt{3}-1)^{-n-2}-(-\sqrt{3}-1)^{-n-2}\right)+ \frac{n!}{\sqrt{3}}\left((\sqrt{3}-1)^{-n-3}-(-\sqrt{3}-1)^{-n-3}\right).$$
\begin{proof} 
It is immediate from the recurrence relation for $L_{n}$ and Theorem \ref{Gmain} that, 
$$L_{n} =n \frac{(n-1)!}{\sqrt{3}}\left((\sqrt{3}-1)^{-n-2}-(-\sqrt{3}-1)^{-n-2}\right) + n(n-1) \frac{(n-2)!}{\sqrt{3}}\left((\sqrt{3}-1)^{-n-3}-(-\sqrt{3}-1)^{-n-3}\right)$$
$$=  \frac{n!}{\sqrt{3}}\left((\sqrt{3}-1)^{-n-2}-(-\sqrt{3}-1)^{-n-2}\right)+ \frac{n!}{\sqrt{3}}\left((\sqrt{3}-1)^{-n-3}-(-\sqrt{3}-1)^{-n-3}\right).$$
\end{proof}
\end{theorem}
We calculate the first 10 terms of the sequence $(L_{n})_{n=0}^{\infty}$ to be, 
$$1,4,15,84,570,4680,44730,488880, 6010200, 82101600.$$

\section{The $\Lambda_{n}^{n}$ case}
Let $n \geq 1$ be a positive integer. In this section we will derive a closed formula for the number of complete $\t$-exceptional sequences in $\module \Lambda_{n}^{n}$. Recall that we denote by $C_n$ the linearly oriented $n$-cycle. 
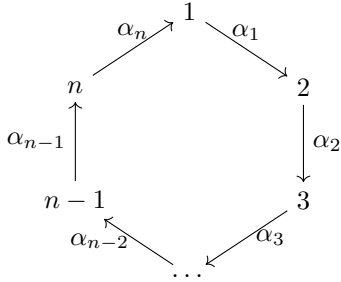
\begin{center}
\begin{tikzpicture}
\node (dots) at (0,0) {$\dots$};
\node (3) at (1.5,1) {$3$};
\node (2) at (1.5,2.5) {$2$};
\node (1) at (0,3.5) {$1$};
\node (n) at (-1.5,2.5) {$n$};
\node (n-1) at (-1.5,1) {$n-1$};

\draw[->] (n-1) -- node[midway, left]{$\alpha_{n-1}$} (n);
\draw[->] (dots) --  node[midway, left]{$\alpha_{n-2}$} (n-1);
\draw[->] (n) --node[midway,above]{$\alpha_{n}$}  (1);
\draw[->] (3) --node[midway,right]{$\alpha_{3}$}  (dots);
\draw[->] (2) --node[midway, right]{$\alpha_{2}$} (3);
\draw[->] (1) --node[midway, above]{$\alpha_{1}$}  (2);
\end{tikzpicture}
\end{center}
The algebra $\Lambda_{n}^{n}$ is defined to be the $\mathbb{F}$-algebra, $\mathbb{F}C_{n}/R_{Q}^{n}$. This is the path algebra of the quiver $C_{n}$ modulo the relations $\alpha_{i}\alpha_{(i+1)_{n}} \dots \alpha_{(i+(n-1))_{n}} = 0$ for $1 \leq j \leq n$.

Let $f_{1} = (2,0)$ and $f_{2} =(1,1)$ in $\mathbb{Z}^{2}$. Recall that Auslander-Reiten quiver of $\module \Lambda_{n}^{n}$ may be identified with the integer lattice, 
$$\mathcal{AR}(\Lambda_{n}^{n}) = \{ af_{1} + bf_{2} : 0 \leq b \leq n-1, \text{ with } a,b \in \mathbb{Z} \}.$$

The following Propositions are needed. 
\begin{proposition}{\label{tauRigid}}{\cite[Proposition 2.5]{Adachi2015tau}}
Let $A$ be a Nakayama algebra. Let $M$ be an indecomposable non-projective module in $\module A$. Then $M$ is
rigid if and only if $l(M)<n$ holds.
\end{proposition}
For our purposes, the following Proposition is a more convenient restatement of Proposition \ref{homAdachi}.

\begin{proposition}{\label{topHom}} Let $M$ be an indecomposable $\Lambda_{n}^{n}$-module with length $1 \leq l(M) \leq n-1$. Then Hom$(X,\t M) \neq 0$ if and only if top$(X) \cong$ top(rad$^{k}(\t M))$ for some $0 \leq k \leq l(M)-1$ and $ l(\text{rad}^{k}(\t M)) \leq l(X).$ 
\begin{proof} All indecomposable modules in $\module \Lambda_{n}^{n}$ have simple tops. 
By Proposition \ref{radAssem}, for a $\Lambda_{n}^{n}$-module $X$, we have that $X = P_{j}/\text{rad}^{l(X)}(P_{j})$ hence top$(X)=S_{j}$. Let $M=P_{i-1} /\text{rad}^{l(M)}P_{i-1}$, then $\t M=P_{(i)_{n}} /\text{rad}^{l(M)}(P_{(i)_{n}})$ by Proposition \ref{radAssem} as well. Observe that for $0 \leq k \leq l(M)-1$,  rad$^{k}(\t M) = P_{(i+k)_{n}}/\text{rad}^{(l(M)-k)}(P_{(i+k)_{n}})$ thus top(rad$^{k}(\t M)) = S_{(i+k)_{n}}$ and $l($rad$^{k}(\t M))=l(M)-k$. By Proposition \ref{homAdachi} we have that, 
$$\text{Hom}(X,\t M) \neq 0 \text{ if and only if } j \in [i , (i+l( M)-1)]_{n}   \text{ and } (i+l( M)-1)_{n} \in [j, (j+l(X)-1)]_{n}.$$

Suppose Hom$(X,\t M) \neq 0$, this implies that $j=(i+k)_{n}$ for some $0 \leq k \leq l(M)-1$, and $(i+l(M)-1)_{n} = (j+a)_{n}$ for some $0 \leq a \leq l(X)-1$. It then immediately follows top$(X) \cong$ top(rad$^{k}(\t M))$ and $ l(\text{rad}^{k}(\t M)) \leq l(X).$ 

For the converse, suppose that top$(X) \cong$ top(rad$^{k}(\t M))$ and $l(\text{rad}^{k}(\t M)) \leq l(X)$. Then $j=(i+k)_{n}$ for some $0 \leq k \leq l(M)-1$. Moreover, $l(\text{rad}^{k}(\t M))=l(M)-k \leq l(X)$ which implies $i+l(M)-1 \leq (i+k)+l(X)-1$ therefore $(i+l( M)-1)_{n} \in [j , (j+l(X)-1)]_{n}$. Hence by Proposition \ref{homAdachi}, Hom$(X,\t M) \neq 0$. 
\end{proof}
\end{proposition}

By Proposition \ref{tauRigid}, every indecomposable module $M$ of $\module \Lambda_{n}^{n}$ is $\t$-rigid in $\module \Lambda_{n}^{n}$ since it is either projective or has length $l(M)<n$. Hence, we once again adopt the same strategy for calculating the number of complete $\t$-exceptional sequences in $\module \Lambda_{n}^{n}$ as we have done thus far. For each $M$ in ind($\module \Lambda_{n}^{n}$), we will calculate the number of complete $\t$-exceptional sequences ending in $M$. By definition a sequence of indecomposable modules $(M_1, M_2, \dots, M_{n-1},M)$ is a $\t$-exceptional sequence in $\module \Lambda_{n}^{n}$ if $M$ is $\t$-rigid and $(M_1, M_2, \dots, M_{n-1})$ is a $\t$-exceptional sequence in $J(M)$. Having seen that every indecomposable $\Lambda_{n}^{n}$-module $M$ is $\t$-rigid, what is left to do is to calculate $J(M)$ for each indecomposable module. Theorem \ref{mainJasso} and Lemma \ref{mainlemma} are once again the main tools these calculations. 

\begin{proposition}{\label{ngammanperp}}Let $M$ be an indecomposable $\Lambda_{n}^{n}$-module with length $1 \leq l(M) \leq n-1$ and top$(M)=S_{i}$. Then for all $1 \leq k \leq l(M)-1$, $$\mathcal{P}(^{\perp}(\t M)) = M \oplus \bigoplus_{s=1}^{l(M)-1} \text{rad}^{s}(M) \oplus \bigoplus_{\substack{1 \leq j \leq n \\  j \notin [i+1,i+l(M)]_{n} }} P_{j}.$$

\begin{proof}
Suppose the $\Lambda_{n}^{n}$-module $M$ has top equal to top$(M)=S_{i}$ and has length $1 \leq l(M) \leq n-1$ i.e. $M$ is not projective. By Proposition \ref{radAssem}, $M = P_{i} / \text{rad}^{l(M)}(P_{i})$ and $\t M = \text{rad}(P_{i})/\text{rad}^{l(M)+1}(P_{i})$ with $l(M)=l(\t M)$. It is easy to see that $\text{top}(\t M) = S_{(i+1)_{n}}$ hence $\t M =  P_{(i+1)_{n}} / \text{rad}^{l(M)}(P_{(i+1)_{n}})$.

By Proposition \ref{topHom}, a $\Lambda_{n}^{n}$-module $X$ is not in  $^{\perp} (\t M)$ if and only if top$(X) \cong$ top(rad$^{k}(\t M))$ for some $0 \leq k \leq l(M)-1$ and $ l(\text{rad}^{k}(\t M)) \leq l(X).$ Let $X = P_{j}/ \text{rad}^{l(X)}(P_{j})$ for some $1 \leq j \leq n$. The statement top$(X) \cong$ top(rad$^{k}(\t M))$ for some $0 \leq k \leq l(M)-1$ means that $j = (i+1+k)_{n}$ for some $0 \leq k \leq l(M)-1$. With this we are able to determine the Ext-projectives in $^{\perp}(\t M)$. 

Let $Y=P_{l}$ be the indecomposable project at the vertex $l$ with $l \neq (i+1+k)_{n}$ for some $0 \leq k \leq l(M)-1$. Then $P_{l}$ is in $ ^{\perp}(\t M)$ by Proposition \ref{topHom}. Moreover Ext$_{\Lambda_{n}^{n}}(P_{l},{^{\perp}(\t M)}) = 0 $ since $P_{l}$ is a projective $\Lambda_{n}^{n}$-module. Hence $P_{l}$ is Ext-projective in $^{\perp}(\t M)$.

Let $Y = \text{rad}^{s}(M)$ for some $1 \leq s \leq l(M)-1$. Then observe that $Y = P_{(i+s)_{n}}/\text{rad}^{(l(M)-s)}(P_{(i+s)_{n}})$ meaning $l(Y)=l(M)-s$. Recall a $\Lambda_{n}^{n}$-module $X$ is not in  $^{\perp} (\t Y)$ if and only if top$(X) \cong$ top(rad$^{r}(\t Y))$ for some $0 \leq r \leq l(M)-s-1$ and $ l(\text{rad}^{r}(\t Y)) \leq l(X).$ Therefore if $X = P_{j}/ \text{rad}^{l(X)}(P_{j})$, then $j= (i+1+s+r)_{n}$ for some $0 \leq r \leq l(M)-s-1$. This implies that $\{ X : \text{Hom}(X, \t Y) \neq 0 \} \subset \{X : \text{Hom}(X, \t M ) \neq 0 \}$, which further implies that Ext$_{\Lambda_{n}^{n}}(Y,N) \cong $ D$\overline{\text{Hom}}_{\Lambda_{n}^{n}}(N, \t Y) = 0$ for all $N$ in $^{\perp}(\t M)$ by the Auslander-Reiten formula. Hence $Y = \text{rad}^{s}(M)$ is an Ext-projective in $^{\perp}(\t M)$.

By Proposition \ref{ExtProj}, $M$ is Ext-projective in $^{\perp}(\t M)$. For every other indecomposable $\Lambda_{n}^{n}$-module $Y$, we have that $\t Y$ is in $^{\perp}(\t M)$, therefore Ext$_{\Lambda_{n}^{n}}(Y, \t Y) \cong \text{D}\overline{\text{Hom}}_{\Lambda_{n}^{n}}(\t Y, \t Y) \neq 0$ i.e. they are not Ext-projective in $^{\perp}(\t M)$. By definition, $T_M = \mathcal{P}(^{\perp}(\t M))$, hence by the above arguments, $$\mathcal{P}(^{\perp}(\t M)) = M \oplus \bigoplus_{s=1}^{l(M)-1} \text{rad}^{s}(M) \oplus \bigoplus_{\substack{\substack{1 \leq j \leq n \\  j \notin [i+1,i+l(M)]_{n}}}} P_{j}.$$
\end{proof}
\end{proposition} 

\begin{proposition}{\label{ngammanProj}}
Let $P_i$ be an indecomposable projective module in $\module \Lambda_{n}^{n}$ for some $1 \leq i \leq n$.  Then the $\t$-perpendicular category of $P_i$ in $\module \Lambda_{n}^{n}$ is $J(P_i) \cong \module \mathbb{A}_{n-1}$, where $\mathbb{A}_{n-1}$ is the Dynkin type $A_{n-1}$ hereditary algebra. 
\begin{proof}
By definition $T_{P_i} = \mathcal{P}(^{\perp}(\t P_i))$. Since $P_i$ is projective, we have that $\t P_i =0$ therefore $^{\perp}(\t P_i) = \module \Lambda_{n}^{n}$. As a result the Ext-projectives of $^{\perp}(\t P_i)$ are just the projectives of $\module \Lambda_{n}^{n}$ therefore $T_{P_i}=\mathcal{P}(^{\perp}(\t P_i)) = \bigoplus_{j=1}^{n}P_{j}$. Thus the $\mathbb{F}$-algebra $E_{P_{i}} = \text{End}_{\Lambda_{n}^{2}}(T_{P_i})$ is precisely given by the path algebra of the quiver $C_{n}^{\text{op}},$
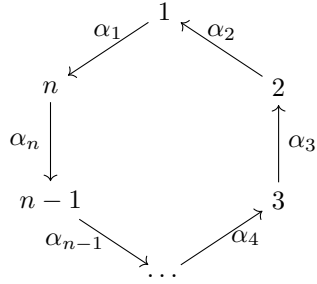
\begin{center}
\begin{tikzpicture}
\node (dots) at (0,0) {$\dots$};
\node (3) at (1.5,1) {$3$};
\node (2) at (1.5,2.5) {$2$};
\node (1) at (0,3.5) {$1$};
\node (n) at (-1.5,2.5) {$n$};
\node (n-1) at (-1.5,1) {$n-1$};

\draw[->] (n) -- node[midway, left]{$\alpha_{n}$} (n-1);
\draw[->] (n-1) --  node[midway, left]{$\alpha_{n-1}$} (dots);
\draw[->] (1) --node[midway,above]{$\alpha_{1}$}  (n);
\draw[->] (dots) --node[midway,right]{$\alpha_{4}$}  (3);
\draw[->] (3) --node[midway, right]{$\alpha_{3}$} (2);
\draw[->] (2) --node[midway, above]{$\alpha_{2}$}  (1);
\end{tikzpicture}
\end{center}

modulo the relations $\alpha_{i}\alpha_{(i+1)_{n}} \dots \alpha_{(i+(n-1))_{n}} =0$ for $1 \leq j \leq n$.

By Lemma \ref{mainlemma}, $D_M = E_{M}/ \langle e_M \rangle$ is the path algebra of the quiver $C_{n}^{\text{op}(i)}$ which is the quiver obtained from $C_{n}^{\text{op}}$ by removing the vertex $i$. More precisely, $C_{n}^{\text{op}(i)}$ is the quiver,

\begin{center}
\begin{tikzcd}
   i+ 1 & i+2 \arrow{l}[swap]{\alpha_{i+2}} & \dots \arrow{l}[swap]{\alpha_{i+3}} & n-1 \arrow{l}[swap]{\alpha_{n-1}} & n \arrow{l}[swap]{\alpha_{n}} & 1 \arrow{l}[swap]{\alpha_{1}} & \dots \arrow{l}[swap]{\alpha_{2}} & i-2 \arrow{l}[swap]{\alpha_{i-2}} & i-1 \arrow{l}[swap]{\alpha_{i-1}}
\end{tikzcd}
\end{center}
with no relations.  It is easy to see that the path algebra $\mathbb{F}C_{n}^{\text{op}(i)}$ is isomorphic to $\mathbb{A}_{n-1}$. Hence the Proposition follows by Theorem \ref{mainJasso}.
\end{proof}
\end{proposition}

\begin{proposition}{\label{ngammanNonProj}}Let $M$ be an indecomposable $\Lambda_{n}^{n}$-module with length $1 \leq l(M) \leq n-1$ and top$(M)=S_{i}$. Which is to say that $M = P_{i}/\text{rad}^{l(M)}(P_{i})$. Then the $\t$-perpendicular category of $M$ in $\module \Lambda_{n}^{n}$ is $J(M)\cong \module \mathbb{A}_{l(M)-1} \oplus \text{mod}\Lambda_{n-l(M)}^{n-l(M)}$, where $\mathbb{A}_{m}$ is the Dynkin type $A_m$ hereditary algebra. 
\begin{proof}
 By Proposition \ref{ngammanperp} the Bongartz completion of $M$ is, $$T_{M} = M \oplus \bigoplus_{s=1}^{l(M)-1} \text{rad}^{s}(M) \oplus \bigoplus_{\substack{1 \leq j \leq n \\  j \notin [i+1,i+l(M)]_{n} }} P_{j}.$$ Hence the $\mathbb{F}$-algebra $E_{M} = \text{End}_{\Lambda_{n}^{n}}(T_M)$ is the path algebra of the following quiver $Q_{n}$.
 
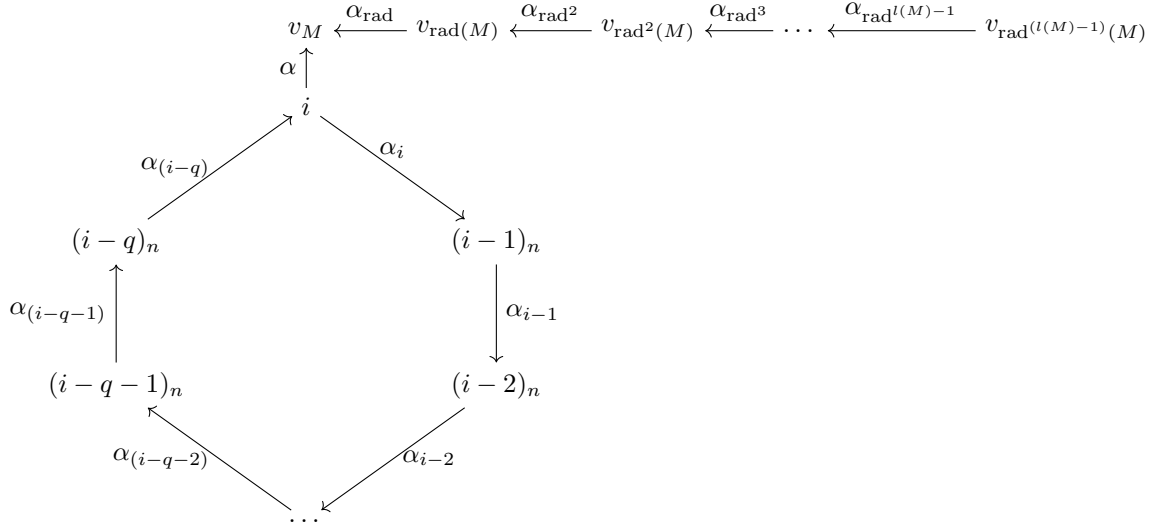
 \begin{center}
\begin{tikzpicture}
\node (dots) at (0,-1) {$\dots$};
\node (i-2) at (2.5,0.8) {$(i-2)_{n}$};
\node (i-1) at (2.5,2.7) {$(i-1)_{n}$};
\node (i) at (0,4.5) {$i$};
\node (n) at (-2.5,2.7) {$(i-q)_{n}$};
\node (n-1) at (-2.5,0.8) {$(i-q-1)_{n}$};
\node(vm) at (0,5.5) {$v_{M}$};
\node(radm) at (2,5.5) {$v_{\text{rad}(M)}$};
\node (rad2m) at (4.5,5.5) {$v_{\text{rad}^{2}(M)}$};
\node (raddots) at (6.5,5.5) {$\dots$};
\node (radlmm) at (10,5.5) {$v_{\text{rad}^{(l(M)-1)}(M)}$};

\draw[->] (n-1) -- node[midway, left]{$\alpha_{(i-q-1)}$} (n);
\draw[->] (dots) -- node[midway, left]{$\alpha_{(i-q-2)}$} (n-1);
\draw[->] (n) --node[midway,left]{$\alpha_{(i-q)}$}  (i);
\draw[->] (i-2) --node[midway,right]{$\alpha_{i-2}$}  (dots);
\draw[->] (i-1) --node[midway, right]{$\alpha_{i-1}$} (i-2);
\draw[->] (i) --node[midway, above]{$\alpha_{i}$}  (i-1);
\draw[->] (i) -- node[midway,left]{$\alpha$} (vm);
\draw[->] (radm) -- node[midway,above]{$\alpha_{\text{rad}}$} (vm);
\draw[->] (rad2m) -- node[midway,above]{$\alpha_{\text{rad}^{2}}$} (radm);
\draw[->] (raddots) -- node[midway,above]{$\alpha_{\text{rad}^{3}}$} (rad2m);
\draw[->] (radlmm) -- node[midway,above]{$\alpha_{\text{rad}^{l(M)-1}}$} (raddots);
\end{tikzpicture}
\end{center}
modulo the relations $\alpha_{r}\alpha_{r-1} \dots \alpha_{r-(n-l(M)-1)} = 0$ for $r \in [i,i-(n-l(M)-1)]_{n}$ and where $q=n-l(M)-1$.

 Let $Q_{n}^{(v_{M})}$ be the quiver obtained from $Q_n$ by removing the vertex $v_{M}$ and any arrows incident to $v_{M}$. More precisely, $Q_{n}^{(v_{M})}$ is the following quiver with two connected components,

 
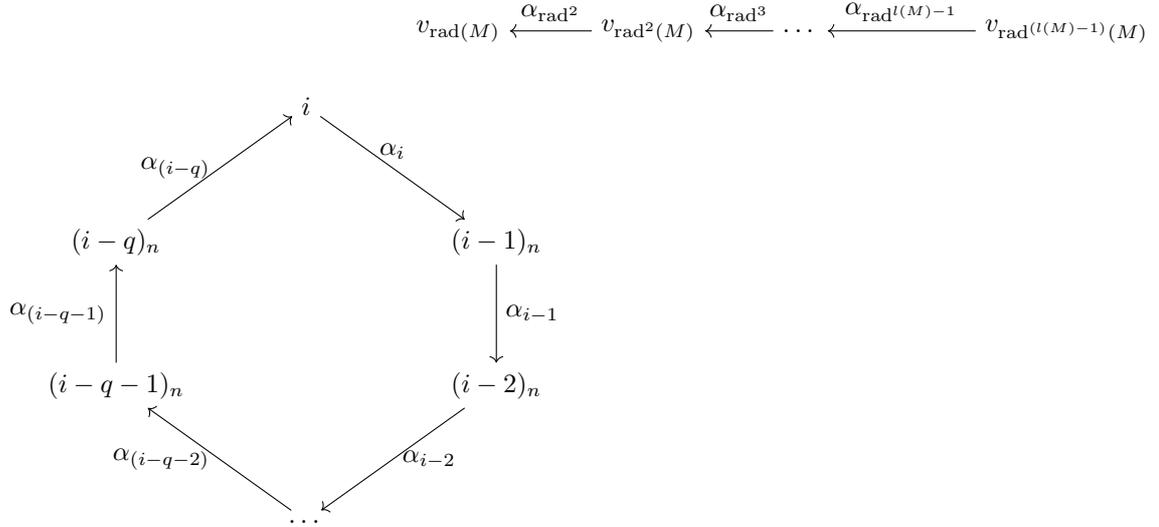
\begin{center}
\begin{tikzpicture}
\node (dots) at (0,-1) {$\dots$};
\node (i-2) at (2.5,0.8) {$(i-2)_{n}$};
\node (i-1) at (2.5,2.7) {$(i-1)_{n}$};
\node (i) at (0,4.5) {$i$};
\node (n) at (-2.5,2.7) {$(i-q)_{n}$};
\node (n-1) at (-2.5,0.8) {$(i-q-1)_{n}$};
\node(radm) at (2,5.5) {$v_{\text{rad}(M)}$};
\node (rad2m) at (4.5,5.5) {$v_{\text{rad}^{2}(M)}$};
\node (raddots) at (6.5,5.5) {$\dots$};
\node (radlmm) at (10,5.5) {$v_{\text{rad}^{(l(M)-1)}(M)}$};

\draw[->] (n-1) -- node[midway, left]{$\alpha_{(i-q-1)}$} (n);
\draw[->] (dots) -- node[midway, left]{$\alpha_{(i-q-2)}$} (n-1);
\draw[->] (n) --node[midway,left]{$\alpha_{(i-q)}$}  (i);
\draw[->] (i-2) --node[midway,right]{$\alpha_{i-2}$}  (dots);
\draw[->] (i-1) --node[midway, right]{$\alpha_{i-1}$} (i-2);
\draw[->] (i) --node[midway, above]{$\alpha_{i}$}  (i-1);
\draw[->] (rad2m) -- node[midway,above]{$\alpha_{\text{rad}^{2}}$} (radm);
\draw[->] (raddots) -- node[midway,above]{$\alpha_{\text{rad}^{3}}$} (rad2m);
\draw[->] (radlmm) -- node[midway,above]{$\alpha_{\text{rad}^{l(M)-1}}$} (raddots);
\end{tikzpicture}
\end{center}
and with relations $\alpha_{r}\alpha_{r-1} \dots \alpha_{r-(n-l(M)-1)} = 0$ for $r \in [i,i-(n-l(M)-1)]_{n}$. By Lemma \ref{mainlemma}, $D_M = E_{M}/ \langle e_M \rangle$ is the path algebra of $Q_{n}^{(v_{M})}$ modulo relations. So it follows that $\module D_{M} \cong \module \mathbb{A}_{l(M)-1} \oplus \text{mod}\Lambda_{n-l(M)}^{n-l(M)}.$ So by Theorem \ref{mainJasso}, the statement of the Proposition follows. 
\end{proof}
\end{proposition}

\begin{theorem}Let $H_{n}$ denote the number of complete $\t$-exceptional sequences in $\module \Lambda_{n}^{n}$. Then $H_{n}$ satisfies the recurrence relation, 
$$ H_{n} =n \sum_{i=1}^{n} \binom{n-1}{i-1} i^{i-2}H_{n-i},$$
with $H_{0} = 1$.

\begin{proof} Let $M$ be an indecomposable $\Lambda_{n}^{n}$-module. Suppose $(M_1, M_2, \dots, M_{n-1}, M)$ is a complete $\t$-exceptional sequence in $\module \Lambda_{n}^{n}$ ending in $M$. Then by definition and the fact that $\delta(J(M))=n-1$, the sequence $(M_1, M_2, \dots, M_{n-1})$ is a complete $\t$-exceptional sequence in $J(M)$. It then follows that the number of complete $\t$-exceptional sequences ending in $M$ is equal to the number of complete $\t$-exceptional sequences in $J(M)$.

The length of $M$ is $1 \leq l(M) \leq n$. For each possible value of $l(M)$, there are $n$ indecomposable $\Lambda_{n}^{n}$-modules of that length. If $l(M)=n$ then $M$ is projective and by Proposition \ref{ngammanProj}, $J(M) \cong  \module \mathbb{{A}}_{n-1}$.  The number of $\t$-exceptional sequences in $\module \mathbb{A}_{n-1}$ was shown in [\cite{seidel2001exceptional} [Proposition 1.1]] to be $n^{n-2}= \binom{n-1}{n-1}n^{n-2}H_{0}$, where $H_{0}=1$. 

If $1 \leq l(M) \leq n-1$, then  by Proposition \ref{ngammanNonProj}, the $\t$-perpendicular category of $M$ is $J(M) \cong \module \mathbb{A}_{l(M)-1} \oplus \module \Lambda_{n-l(M)}^{n-l(M)}$. Arguing as above, the number of complete $\t$-exceptional sequences ending in $M$ is equal to the number of complete $\t$-exceptional sequences in $ \module \mathbb{A}_{l(M)-1} \oplus \module \Lambda_{n-l(M)}^{n-l(M)}$. By Theorems \ref{main1} and \ref{main2}, this is equal to $$\binom{n-1}{n-l(M) , l(M)-1}l(M)^{(l(M)-2)}H_{n-l(M)} = \binom{n-1}{l(M)-1}l(M)^{(l(M)-2)}H_{n-l(M)}.$$
So it follows that, $$H_{n}= \sum_{l(M)=1}^{n} n \binom{n-1}{l(M)-1} l(M)^{l(M)-2}H_{n-l(M)} =n \sum_{i=1}^{n} \binom{n-1}{i-1} i^{i-2}H_{n-i}.$$ 

It is trivial to see that $H_{1}=1$. Using the recurrence we obtain $H_{1} = \binom{0}{0}1^{-1}H_{0}=1$, therefore $H_{0}=1$. 
\end{proof}
\end{theorem}

We are now in a position to derive the exponential generating function of $H_{n}$. First we state the following results which will be useful in deriving the exponential generating function. 

\begin{lemma}{\label{convolutions}}{\cite[Section 2.3 Rule $3^{\prime}$]{wilf2005generatingfunctionology}} Let $f = \sum_{n=0}^{\infty} a_{n}\frac{x^{n}}{n!}$ and $g = \sum_{n=0}^{\infty} b_{n}\frac{x^{n}}{n!}$ be the generating functions of the sequences $\{ a_{n} \}_{n=0}^{\infty}$ and $\{ b_{n} \}_{n=0}^{\infty}$ respectively. Then the series $fg$ is the exponential generating function of the series, 
$$ \Bigg \{ \sum_{k=0}^{n} \binom{n}{k}a_{k}b_{n-k} \Bigg\}_{n=0}^{\infty}.$$
\end{lemma}
\begin{definition}{\cite{DWalsh}}
Let $f : \{1,2, \dots, n\} \rightarrow \{1,2, \dots, n+m\}$ be a function. For a positive integer $1 \leq k \leq n$, a $k$-cycle of $f$ is a sequence of distinct elements $(x_1, x_2, \dots, x_k)$ in the domain of $f$, such that $f(x_{i})=x_{(i+1)_{n}}$. The map $f$ is called an \textit{acyclic function} if it does not have any $k$-cycles for all $1 \leq k \leq n$. 
\end{definition}

\begin{lemma}{\label{LambertW}}{\cite[Theorem 1]{DWalsh}} Let $a \geq 1$ and $b \geq 0$ be integers. The number $N(a,b)$ of acyclic functions from domain $D= \{1,2,\dots,b \}$ to codomain $C = \{1,2, \dots ,a+b\} $ is given by,
$$N(a,b) = a(a+b)^{b-1}.$$
Moreover, for a fixed positive integer $a$, the exponential generating function of $N(a,b)$ is given by,
$$\sum_{b=0}^{\infty} N(a,b)\frac{x^{b}}{b!} = e^{-aW(-x)},$$
where $W$ is Lambert's $W$ function. 
\end{lemma}
The Lambert $W$ function is defined to be the inverse of the function where $w \mapsto we^{w}$. The function $W$ has many applications in mathematics. For example, it is used in the enumeration of trees and the calculation of water-wave heights. The reader is referred to \cite{corless1996lambertw} for more on Lambert's $W$ function. 

\begin{theorem}The exponential generating function of $H_{n}$ is,
$$ \sum_{n=0}^{\infty} H_{n}\frac{x^{n}}{n!} = \frac{1}{1+W(-x)}, $$
where $W$ is Lambert's $W$ function and $H_{n}$ is given by the closed formula, $$H_{n} = n^{n}.$$
\begin{proof}
Let $a_{n}$ be the sequence $a_{n} = (n+1)^{n-1}$. Let $h(x) = \sum_{n=0}^{\infty} H_{n}\frac{x^{n}}{n!}$ and $g(x) = \sum_{n=0}^{\infty} a_{n}\frac{x^{n}}{n!}$ be exponential generating functions of $H_{n}$ and $a_{n}$ respectively. 
Recall the recurrence relation of $H_{n}$ is given by, 
$$ H_{n} =n \sum_{k=1}^{n} \binom{n-1}{k-1} k^{k-2}H_{n-k},$$
so, $$\frac{H_{n}}{n} = \sum_{k=1}^{n} \binom{n-1}{k-1} k^{k-2}H_{n-k}.$$
We make the change of variable $j=k-1$ in $\frac{H_{n}}{n}$ to obtain the following. 
$$ \frac{H_{n}}{n} = \sum_{j=0}^{n-1} \binom{n-1}{j} (j+1)^{j-1}H_{n-(j+1)},$$
thus 
$$ \frac{H_{n+1}}{n+1} = \sum_{j=0}^{n} \binom{n}{j} (j+1)^{j-1}H_{n-j}.$$
We now study the exponential generation function of $\frac{H_{n+1}}{n+1}$,
$$\sum_{n=0}^{\infty} \frac{H_{n+1}}{n+1}\frac{x^{n}}{n!} = \sum_{n=0}^{\infty} \left( \sum_{j=0}^{n} \binom{n}{j} (j+1)^{j-1}H_{n-j}\right)\frac{x^{n}}{n!}.$$
By Lemma \ref{convolutions}, the right hand side is given the product $g(x)h(x)$. So we have, 
$$\sum_{n=0}^{\infty} \frac{H_{n+1}}{n+1}\frac{x^{n}}{n!} = g(x)h(x).$$
We can manipulate the right hand side so that the exponent of $x$ matches the factorial, hence 
$$\frac{1}{x} \sum_{n=0}^{\infty} H_{n+1}\frac{x^{n+1}}{(n+1)!} =g(x)h(x),$$
so we can write the left hand side in terms of $h(x)$ as follows, 
$$\frac{1}{x}(h(x)-H_{0})=g(x)h(x).$$
Since $H_{0}=1$, 
$$h(x)-1 = xh(x)g(x).$$
By Lemma \ref{LambertW}, $g(x)= e^{-W(-x)}$ therefore,
 $$h(x) = \frac{1}{1-xg(x)} = \frac{1}{1 -xe^{-W(-x)}}.$$
 Recall that Lambert's $W$ function is defined by the equation $x=W(x)e^{W(x)}$ (See \cite{corless1996lambertw} for more on Lambert's $W$ function), thus $-xe^{-W(-x)}=W(-x)$, giving us that, 
 $$h(x) = \frac{1}{1+W(-x)} = \frac{1}{1-T(x)},$$
 where $T(x)=-W(-x)$ is called Euler's tree function, again see \cite{corless1996lambertw}. 
This exponential generating function is precisely the exponential generating function of the sequence $n^n$, see; \cite{knuth1989recurrence} Section 2 equation 2.7 and \cite{riordan1968combinatorial}.
\end{proof}
\end{theorem}
It is interesting to note that $n^n$ also counts the number of  complete exceptional sequences the Dykin algebras of type B and C; see section 5 of \cite{obaid2013number}.

\section{The $\Gamma_{n}^{n-1}$ case}
Let $n \geq 1$ be a positive integer. In this section we will study the combinatorics for the number of complete $\t$-exceptional sequence in $\module \Gamma_{n}^{n-1}$. Recall that we denote by $A_n$ the linearly oriented quiver with $n$ vertices, 
\begin{center}
\begin{tikzcd}
    1 \arrow[r,"\alpha_{1}"]& 2 \arrow[r,"\alpha_{2}"] & 3 \arrow[r, "\alpha_{3}"] & \dots \arrow[r,"\alpha_{n-2}"] & n-1 \arrow[r,"\alpha_{n-1}"] & n
\end{tikzcd}.
\end{center} The algebra $\Gamma_{n}^{n-1}$ is defined to be the $\mathbb{F}$-algebra, $\mathbb{F}A_{n}/R_{Q}^{n-1}$. This is the path algebra of the quiver $A_{n}$ modulo the relation $\alpha_{1}\alpha_{2}\dots\alpha_{n-1}=0$.

Let $f_{1} = (2,0)$ and $f_{2} =(1,1)$ in $\mathbb{Z}^{2}$. The Auslander-Reiten quiver of $\module \Gamma_{n}^{n-1}$ may be identified with the integer lattice, 
$$\mathcal{AR}(\Gamma_{n}^{n-1}) = \{  af_{1} + bf_{2} : 0 \leq a \leq n-1, 0 \leq b \leq n-2 \text{ such that } a,b \in \mathbb{Z} \text{ and } a+b \leq n-1  \}.$$

Observe the following. Let $M$ be an indecomposable module in $\module \Gamma_{n}^{n-1}$, then $M$ belongs to one of the following disjoint sets. The first set contains the indecomposable projective modules $P_{j}$ for $1 \leq j \leq n$. The second set contains non-projective modules of the form $M = $rad$^{i}(P_{1})$ where $1 \leq i \leq n-2$ and $P_{1}$ is the indecomposable projective at vertex 1. The third set contains indecomposable modules which are neither projective or of the form $M = $rad$^{i}(P_{1})$ for $1 \leq i \leq n-2$. Any indecomposable module $M$ in $\module \Gamma_{n}^{n-1}$ has length $l(M) < n$, therefore by Proposition \ref{tauRigid}, every indecomposable module of $\module \Gamma_{n}^{n-1}$ is $\t$-rigid. 

\begin{proposition}{\label{gamma(n-1)proj}}
Let $P_{i}$ be an indecomposable projective module in $\module \Gamma_{n}^{n-1}$ for some $1 \leq i \leq n$. Then the $\t$-perpendicular category of $P_{i}$ in $\module \Gamma_{n}^{n-1}$ is $J(P_{i}) \cong \module \mathbb{A}_{n-i} \oplus \module \mathbb{A}_{i-1},$ where $\mathbb{A}_{j}$ is the Dynkin type $A_j$ hereditary algebra. 
\begin{proof} Let $P_{i}$ be an indecomposable projective with length $1 \leq l(P_{i}) \leq n-1$. By definition the Bongartz completion $T_{P_{i}} = \mathcal{P}({^{\perp}(\t P_{i})}).$ Since $P_{i}$ is projective, $\t P_{i} = 0$ therefore $^{\perp}(\t P_{i}) = \module \Gamma_{n}^{n-1}$, hence the Bongartz completion $T_{P_{i}} = \bigoplus_{j=1}^{n} P_{j}$. Thus the $\mathbb{F}$-algebra $E_{P_{i}} = \text{End}_{\Gamma_{n}^{n-1}}(T_{P_{i}})$ is precisely the algebra $\Gamma_{n}^{n-1}$, the path algebra of the quiver $A_{n}^{\text{op}}$, 
\begin{center}
\begin{tikzcd}
    1 & 2 \arrow{l}[swap]{\alpha_{2}} & 3 \arrow{l}[swap]{\alpha_{3}} & \dots \arrow{l}[swap]{\alpha_{4}} & i-1 \arrow{l}[swap]{\alpha_{i-1}} & i \arrow{l}[swap]{\alpha_{i}} & i+1 \arrow{l}[swap]{\alpha_{i+1}} & \dots \arrow{l}[swap]{\alpha_{i+2}} & n-1 \arrow{l}[swap]{\alpha_{n-1}} & n \arrow{l}[swap]{\alpha_{n}}
\end{tikzcd}
\end{center}
modulo the relation $\alpha_{n}\alpha_{n-1}\dots\alpha_{1}=0$.
Let $A_{n}^{\text{op}(i)}$ be the quiver obtained from $A_{n}^{\text{op}}$ by removing the vertex $i$ and all arrows incident to $i$.
\begin{center}
\begin{tikzcd}
    1 & 2 \arrow{l}[swap]{\alpha_{2}} & 3 \arrow{l}[swap]{\alpha_{3}} & \dots \arrow{l}[swap]{\alpha_{4}} & i-1 \arrow{l}[swap]{\alpha_{i-1}} & & i+1  & \dots \arrow{l}[swap]{\alpha_{i+2}} & n-1 \arrow{l}[swap]{\alpha_{n-1}} & n \arrow{l}[swap]{\alpha_{n}}
\end{tikzcd}
\end{center}
The quiver $A_{n}^{\text{op}(i)}$ has no relations. By Lemma \ref{mainlemma}, $D_M = E_{M}/ \langle e_M \rangle$ is the path algebra of the quiver $A_{n}^{\text{op}(i)}$. Since $D_M = \mathbb{F}A_{n}^{\text{op}(i)}$, it follows that $J(M) \cong \module \mathbb{A}_{n-i} \oplus \module \mathbb{A}_{i-1}$ by Theorem \ref{mainJasso}.
\end{proof}
\end{proposition}

\begin{proposition}{\label{gamma(n-1)rad}}
Let $M$ be an indecomposable module in $\module \Gamma_{n}^{n-1}$ of the form $M = \text{rad}^{i}(P_{1})$ for some $ 1 \leq i \leq n-2$ with length $1 \leq l(M) \leq n-2.$ Then the $\t$-perpendicular category of $M$ in $\module \Gamma_{n}^{n-1}$ is
\[ J(M) \cong \begin{cases} 
           \module \mathbb{A}_{l(M)-1} \oplus \module \mathbb{A}_{1} \oplus \module \mathbb{A}_{1} & i=1 \\
          \module \mathbb{A}_{n-l(M)} \oplus \module \mathbb{A}_{l(M)-1} & i\neq 1 \\
       \end{cases},
    \]
 where $\mathbb{A}_{j}$ is the Dynkin type $A_j$ hereditary algebra.
\begin{proof} Consider $P_{j}$ the indecomposable projective module at the vertex $j$ in $\module \Gamma_{n}^{n-1}$ with $j \neq 1$. Then it is easy to see that $P_{j}=\text{rad}^{k}(P_2)$ for some $0 \leq k\leq n-2$ and that $P_{j}$ has length $l(P_{j})=n-j+1$. From this it follows that rad$^{q}(P_{j}) = P_{j+q}$ where $0 \leq q \leq n-j$.

Let $M = \text{rad}^{i}(P_{1})$ for some $1 \leq i \leq n-1$. We observe that $l(M) = n-i-1$. In accordance to Proposition \ref{radAssem}, $M$ may in fact be written as $M = P_{i+1} / \text{rad}^{n-i-1}(P_{i+1})$. Using Proposition \ref{radAssem} again, we can see that Auslander-Reiten translate of $M$ is given by $\t M = \text{rad}(P_{i+1})/\text{rad}^{n-i}(P_{i+1})=P_{i+2}$ because $\text{rad}(P_{i+1})=P_{i+2}$ and $l(P_{i+1})=n-(i+1)+1=n-i$, hence rad$^{n-i}(P_{i+1})=0$. So we see that the only indecomposable $\Gamma_{n}^{n-1}$-modules not in $^{\perp}( \t M)$ are the projectives $P_{j} = \text{rad}^{s}(P_{i+2})$ for $0 \leq s < n-i-1$, in other words $i+2 \leq j \leq n$ since $l(P_{i+2})=n-i-1$.

We are now in the position to determine the Ext-projectives of $^{\perp} (\t M) \text{ where } M = \text{rad}^{i}(P_{1})$. By the above calculation, we can say that for $1 \leq j \leq i+1$, the projective $P_{j} \text{  is in } ^{\perp} (\t M)$ hence Ext$_{\Gamma_{n}^{n-1}}( P_{j}, ^{\perp} (\t M))=0$. 

Let $N = \text{rad}^{j}(P_{1})$ for some $j > i$. Arguing as above we can see that the only indecomposable $\Gamma_{n}^{n-1}$-modules not in $^{\perp}( \t N)$ are the indecomposable projectives $P_{m}$ where $j+2 \leq m \leq n$, so it follows that $\{X: \text{Hom}(X, \t N) \neq 0 \} \subset \{ X : \text{Hom}(X, \t M ) \neq 0 \}$. This implies that Ext$_{\Gamma_{n}^{n-1}}(N,X) \cong $ D$\overline{\text{Hom}}_{\Gamma_{n}^{n-1}}(X, \t N) = 0$ for all $X$ in $^{\perp}(\t M)$ by the Auslander-Reiten formula. Hence $N = \text{rad}^{j}(P_{1})$ is an Ext-projective in $^{\perp}(\t M)$.

For every other indecomposable $\Gamma_{n}^{n-1}$ module $Y$, we have that $\t Y$ is in $^{\perp}(\t M)$, therefore since Ext$_{\Gamma_{n}^{n-1}}(Y, \t Y) \cong \text{D}\overline{\text{Hom}}_{\Gamma_{n}^{n-1}}(X, \t N)  \neq 0$. Therefore these modules are not Ext-projective in $^{\perp}(\t M)$. By definition $T_M = \mathcal{P}(^{\perp}(\t M))$, so by the above arguments, $$\mathcal{P}(^{\perp}(\t M)) = M \oplus \bigoplus_{s=i+1}^{n-2} \text{rad}^{s}(P_{1}) \oplus \bigoplus_{j=1}^{i+1} P_{j}.$$
In the case when $i \neq 1$ the $\mathbb{F}$-algebra $E_{M} = \text{End}_{\Gamma_{n}^{n-1}}(T_{M})$ is the path algebra of the quiver $Q_{n}$,
\begin{center}
\begin{tikzcd}[column sep=3.5em,row sep=3.5em]
 &  &  & & i+1 \arrow[r,"\alpha_{i+1}"] \arrow[d, "\alpha"] & i \arrow[r,"\alpha_{i}"] & \dots \arrow[r,"\alpha_{2}"] & 1 \\ 
 v_{\text{rad}^{n-2}} \arrow[r,"\alpha_{\text{rad}^{n-2}}"] &v_{\text{rad}^{n-3}} \arrow[r,"\alpha_{\text{rad}^{n-3}}"]  & \dots \arrow[r,"\alpha_{\text{rad}^{i+1}}"] & v_{\text{rad}^{i+1}} \arrow[r,"\alpha_{\text{rad}^{i+1}}"] & v_{M} \arrow[bend right =20, swap]{urrr}{\alpha_{V_{M}}} & & 
\end{tikzcd}.
\end{center}
By Lemma \ref{mainlemma}, $D_M = E_{M}/ \langle e_M \rangle$ is the path algebra of the quiver $Q_{n}^{(v_M)}$ which is the quiver obtained from $Q_{n}$ by removing the vertex $v_{M}$ and all arrows incident to $v_{M}$. The quiver $Q_{n}^{(v_{M})}$ has two connected components.
\begin{center}
\begin{tikzcd}[column sep=3.5em,row sep=3.5em]
 &  &  & & i+1 \arrow[r,"\alpha_{i+1}"]  & i \arrow[r,"\alpha_{i}"] & \dots \arrow[r,"\alpha_{2}"] & 1 \\ 
 v_{\text{rad}^{n-2}} \arrow[r,"\alpha_{\text{rad}^{n-2}}"] &v_{\text{rad}^{n-3}} \arrow[r,"\alpha_{\text{rad}^{n-3}}"]  & \dots \arrow[r,"\alpha_{\text{rad}^{i+1}}"] & v_{\text{rad}^{i+1}} & & & 
\end{tikzcd}
\end{center}
Since $D_M = \mathbb{F}Q_{n}^{(v_M)}$, it follows that $J(M) \cong \module \mathbb{A}_{i+1} \oplus \module \mathbb{A}_{n-i-2}$ by Theorem \ref{mainJasso}. Recall that $l(M) = n-i-1$, hence $J(M) \cong \module \mathbb{A}_{n-l(M)} \oplus \module \mathbb{A}_{l(M)-1}$.

When $i = 1$ however, the $\mathbb{F}$-algebra $E_{M} = \text{End}_{\Gamma_{n}^{n-1}}(T_{M})$ is the path algebra of the quiver $Q_{n}^{\prime}$,
\begin{center}
\begin{tikzcd}[column sep=3.5em,row sep=3.5em]
 &  &  & & 2 \arrow[d,"\alpha"] & & &  \\ 
 v_{\text{rad}^{n-2}} \arrow[r,"\alpha_{\text{rad}^{n-2}}"] &v_{\text{rad}^{n-3}} \arrow[r,"\alpha_{\text{rad}^{n-3}}"]  & \dots \arrow[r,"\alpha_{\text{rad}^{i+1}}"] & v_{\text{rad}^{i+1}} \arrow[r,"\alpha_{\text{rad}^{i+1}}"] & v_{M} \arrow[r, "\alpha_{v_{M}}"] & 1 & 
\end{tikzcd}
\end{center}
with no relations. By Lemma, \ref{mainlemma} $D_M = E_{M}/ \langle e_M \rangle$ is the path algebra of the quiver $Q_{n}^{{\prime}(v_M)}$ which is the quiver obtained from $Q_{n}^{\prime}$ by removing the vertex $v_{M}$ and all arrows incident to $v_{M}$. This quiver has three connected components.
\begin{center}
\begin{tikzcd}[column sep=3.5em,row sep=3.5em]
 &  &  & & 2 & & &  \\ 
 v_{\text{rad}^{n-2}} \arrow[r,"\alpha_{\text{rad}^{n-2}}"] &v_{\text{rad}^{n-3}} \arrow[r,"\alpha_{\text{rad}^{n-3}}"]  & \dots \arrow[r,"\alpha_{\text{rad}^{i+1}}"] & v_{\text{rad}^{i+1}} &  & 1 & 
\end{tikzcd}.
\end{center}
Since $D_M = \mathbb{F}Q_{n}^{\prime(v_M)}$, it follows that $J(M) \cong \module \mathbb{A}_{n-3} \oplus \module \mathbb{A}_{1} \oplus \module \mathbb{A}_{1}$ by Theorem \ref{mainJasso}. Since $l(\text{rad}^{1}(P_{1})) =n-2$ then $J(M) \cong \module \mathbb{A}_{l(M)-1} \oplus \module \mathbb{A}_{1} \oplus \module \mathbb{A}_{1}$.
\end{proof}
\end{proposition}

\begin{proposition}{\label{gamma(n-1)nrad}}
Let $M$ be an indecomposable $\Gamma_{n}^{n-1}$-module such that $M \neq \text{rad}^{k}(P)$ for some indecomposable projective $P$ and positive integer $k$. Suppose $M$ has length $ 1 \leq l(M) \leq n-2$, then the $\t$-perpendicular category of $M$ in $\module \Gamma_{n}^{n-1}$ is $J(M) \cong \module \mathbb{A}_{l(M)-1} \oplus \module \Gamma_{n-l(M)}^{n-l(M)-1}$.
\begin{proof}
By Proposition \ref{radAssem}, we can write $M = P_{i}/\text{rad}^{l(M)}(P_{i})$ for some $1 \leq i \leq n-1$ and $\t M = \text{rad}(P_{i})/ \text{rad}^{l(M)+1}(P_{i})$. Since rad$(P_{i}) = P_{i+1}$, we have that $\t M = P_{i+1} / \text{rad}^{l(M)}(P_{i+1}).$ Now let $X = P_{j}/\text{rad}^{l}(P_{j})$ be an arbitrary indecomposable $\Gamma_{n}^{n-1}$ module. By Proposition \ref{homAdachi}, Hom$(X, \t M) \neq 0$ if and only $ j \in [i+1,i+l(M)]_{n}$ and $i +l(M) \in [j, j+l-1]_{n}$. From this it follows that $P_{j}$ is not in $^{\perp}( \t M)$ if $i+1 \leq j \leq i+l(M)$. Hence Ext$_{\Gamma_{n}^{n-1}}(P_{j}, ^{\perp}( \t M))=0$ if $j \notin [i+1,i+l(M)]$.

Consider the module $\text{rad}^{s}(M)$ for $1 \leq s \leq l(M)-1$. The length of $\text{rad}^{s}(M)$ is given by $l(\text{rad}^{s}(M)) = l(M)-s$. Moreover, rad$^{s}(M) = P_{i+s}/ \text{rad}^{l(M)-s}(P_{i+s})$, from which it follows that $\t \text{rad}^{s}(M) = P_{i+s+1} / \text{rad}^{l(M)-s}(P_{i+s+1}).$ Again let $X = P_{j}/\text{rad}^{l}(P_{j})$ be an arbitrary indecomposable $\Gamma_{n}^{n-1}$ module. By Proposition \ref{homAdachi}, Hom$(X, \t\ \text{rad}^{s}(M) ) \neq 0$ if and only $ j \in [i+s+1,i+l(M)]_{n}$ and $i +l(M) \in [j, j+l-1]_{n}$. Therefore $\{ X : \text{Hom}(X, \t \text{rad}^{s}(M)) \neq 0 \} \subset \{ X: \text{Hom}(X, \t M) \neq 0 \}$, which implies that Ext$_{\Gamma_{n}^{n-1}}(\text{rad}^{s}(M), Y) \cong $ D$\overline{\text{Hom}}_{\Gamma_{n}^{n-1}}(Y, \t \text{rad}^{s}(M)) = 0$ for all $Y$ in $^{\perp}(\t M)$. In other words, rad$^{s}(M)$ is Ext-projective in $^{\perp}(\t M)$.

By Proposition \ref{ExtProj}, $M$ is Ext-projective in $^{\perp}(\t M)$, so $$\mathcal{P}(^{\perp}( \t M)) = M \oplus \bigoplus_{s=1}^{l(M)-1} \text{rad}^{s}(M) \oplus \bigoplus_{j \notin[ i+1,i+l(M)]} P_{j}.$$ By definition, the Bongartz completion $T_{M} = \mathcal{P}(^{\perp}( \t M))$, so the $\mathbb{F}$-algebra $E_{M} =\text{End}_{\Gamma_{n}^{n-1}}(T_{M})$ is the path algebra of the quiver $Q_{n}$ modulo relations (set $l(M):=m$),

\begin{center}
\begin{tikzcd}[column sep=3.5em,row sep=3.5em]
 n \arrow[r,"\alpha_{n}"]& n-1 \arrow[r,"\alpha_{n-1}"]  & \dots \arrow[r, "\alpha_{i+m+2}"] & i+m+1 \arrow[r,"\alpha_{i+m+1}"] & i \arrow[r,"\alpha_{i}"] \arrow[d, "\alpha"] & i-1 \arrow[r,"\alpha_{i-1}"] & \dots \arrow[r,"\alpha_{2}"] & 1 \\ 
 v_{\text{rad}^{m-1}} \arrow[r,"\alpha_{\text{rad}^{m-1}}"] &v_{\text{rad}^{m-2}} \arrow[r,"\alpha_{\text{rad}^{m-2}}"]  & \dots \arrow[r,"\alpha_{\text{rad}^2}"] & v_{\text{rad}^{1}} \arrow[r,"\alpha_{\text{rad}^{1}}"] & v_{M} & &
\end{tikzcd}.
\end{center}
Since the vertices of the top row of the quiver correspond to the indecomposable projectives of $\module \Gamma_{n}^{n-1}$ and the arrows reflect the relations the corresponding maps between the projectives, we see that we have the relation $\alpha_{n}\alpha_{n-1} \dots \alpha_{i+m+1}\alpha_{i}\alpha_{i-1}\dots \alpha_{1}=0$. Let $Q_{n}^{(v_{M})}$ be the quiver obtained from $Q_{n}$ by removing the vertex $v_{M}$ and all the arrows incident to $v_{M}$,

\begin{center}
\begin{tikzcd}[column sep=3.5em,row sep=3.5em]
 n \arrow[r,"\alpha_{n}"]& n-1 \arrow[r,"\alpha_{n-1}"]  & \dots \arrow[r, "\alpha_{i+m+2}"] & i+m+1 \arrow[r,"\alpha_{i+m+1}"] & i \arrow[r,"\alpha_{i}"]& i-1 \arrow[r,"\alpha_{i-1}"] & \dots \arrow[r,"\alpha_{2}"] & 1 \\ 
 v_{\text{rad}^{m-1}} \arrow[r,"\alpha_{\text{rad}^{m-1}}"] &v_{\text{rad}^{m-2}} \arrow[r,"\alpha_{\text{rad}^{m-2}}"]  & \dots \arrow[r,"\alpha_{\text{rad}^2}"] & v_{\text{rad}^{1}} & & 
\end{tikzcd}
\end{center}
with the relation $\alpha_{n}\alpha_{n-1} \dots \alpha_{i+m+1}\alpha_{i}\alpha_{i-1}\dots \alpha_{1}=0$. By Lemma \ref{mainlemma}, $D_M = E_{M}/ \langle e_M \rangle$ is the path algebra of the quiver $Q_{n}^{(v_{M})}$ modulo the relation $\alpha_{n}\alpha_{n-1} \dots \alpha_{i+m+1}\alpha_{i}\alpha_{i-1}\dots \alpha_{1}=0$. It follows that $J(M) \cong \module \mathbb{A}_{l(M)-1} \oplus \module \Gamma_{n-l(M)}^{n-l(M)-1}$ by Theorem \ref{mainJasso}. 

\end{proof}
\end{proposition}

\begin{theorem}{\label{kn-reccurence}}Let $K_{n}$ denote the number of complete $\t$-exceptional sequences $\module \Gamma_{n}^{n-1}$. Then $K_{n}$ satisfies the recurrence relation;
$$ K_{n} =\sum_{i=1}^{n} \binom{n-1}{i-1}(n-i+1)^{(n-i-1)} i^{i-2} + \sum_{i=1}^{n-3} \binom{n-1}{i-1}(n-i+1)^{(n-i-1)} i^{i-2} + (n-1)(n-2)^{(n-3)}$$ $$ + \sum_{i=1}^{n-2} \binom{n-1}{i-1}(n-i-1) i^{i-2}K_{n-i}$$
with $K_{1} = 1$.
\begin{proof}
Let $M$ be an indecomposable module in $\module \Gamma_{n}^{n-1}$. Suppose $(X_1, X_2, \dots , X_{n-1},M)$ is a $\t$-exceptional sequence in $\module \Gamma_{n}^{n-1}$. Then by definition and the fact that $\delta(J(M))=n-1$, the sequence $(X_1, X_2, \dots , X_{n-1})$ is a complete $\t$-exceptional sequence in $J(M)$. Hence the number of complete $\t$-exceptional sequences ending in $M$ is equal to the number of complete $\t$-exceptional sequences in $J(M)$. 

Suppose $M$ is projective, hence $M=P_{i}$ for some $1 \leq i \leq n$, then by Proposition \ref{gamma(n-1)proj} the $\t$-perpendicular category $J(M) \cong \module \mathbb{A}_{n-i} \oplus \module \mathbb{A}_{i-1}$. The number of complete $\t$-exceptional sequences in $\module \mathbb{A}_{l}$ is precisely the number of complete exceptional sequence in $\module \mathbb{A}_{l}$ which is shown in [\cite{seidel2001exceptional} [Proposition 1.1]]  to be $(l+1)^{(l-1)}$. Therefore by Theorems \ref{main1} and \ref{main2} the number of complete $\t$-exceptional sequence ending in $M=P_{i}$ is $\binom{n-1}{i-1}(n-i+1)^{n-i-1} i^{i-2}$.

Suppose $M = \text{rad}^{i}(P_{1})$ for some $1 \leq i \leq n-2$. If $i=1$ then we saw in Proposition \ref{gamma(n-1)rad} that $J(M) \cong \module \mathbb{A}_{n-3} \oplus \module \mathbb{A}_{1} \oplus \module \mathbb{A}_{1}$. Arguing as above it follows that the number of complete $\t$-exceptional sequences ending in rad$^{1}(P_{1})$ is $\binom{n-1}{n-3,1,1}(n-2)^{(n-4)}2^{0}2^{0} = (n-1)(n-2)^{(n-3)}.$ If it is the case that $2 \leq i \leq n-2$, then $J(M) \cong \module \mathbb{A}_{n-l(M)} \oplus \module \mathbb{A}_{l(M)-1}$. Therefore the number of complete $\t$-exceptional sequences ending in $M = \text{rad}^{i}(P_{1})$ for some $ 2 \leq i \leq n-2$ is $\binom{n-1}{l(M)-1}(n-l(M)+1)^{n-l(M)-1} l(M)^{l(M)-2}$, where $l(M)$ is the length of $M$.

Finally suppose that $M$ is not of the form rad$^{i}(P)$ for some indecomposable projective module $P$. By Proposition \ref{gamma(n-1)nrad}, $J(M) \cong \module \mathbb{A}_{l(M)-1} \oplus \module \Gamma_{n-l(M)}^{n-l(M)-1}$ where $l(M)$ is the length of $M$. Therefore the number of complete $\t$-exceptional sequences ending in $M$ is $\binom{n-1}{l(M)-1}K_{n-l(M)} l(M)^{l(M)-2}$. Observe that in this case the length of $M$ is $1 \leq l(M) \leq n-2$ and for each fixed value of $l(M)$ there are $n-l(M)-1$ indecomposable modules $M$ such that $M \neq \text{rad}^{i}(P)$.

By counting the number of complete $\t$-exceptional sequences ending in each indecomposable $\Gamma_{n}^{n-1}$-module $M$, the recurrence relation of $K_{n}$ follows. It is also trivial to see that $K_{1}=1$.
\end{proof}
\end{theorem}

\begin{theorem} Let $h(x)= \sum_{n=0}^{\infty} K_{n}\frac{x^{n}}{n!}$ be the exponential generating function of $K_{n}$. 
Then $h(x)$ satisfies the first order linear ODE, 
$$h^{\prime}(x)(1-xe^{-W(-x)}) + h(x)e^{-W(-x)}= 2e^{-2W(-x)} - e^{-W(-x)} + W(-x) + \frac{1}{2}xW(-x).$$
\begin{proof}
Let $h(x) = \sum_{n=0}^{\infty} K_{n}\frac{x^{n}}{n!} $ be the exponential generating function of $K_{n}$. Let $a(n) = (n+1)^{n-1}$. Let $g(x) = \sum_{n=0}^{\infty} (n+1)^{n-1}\frac{x^{n}}{n!}$. Then $g(x) = e^{-W(-x)}$ by Lemma \ref{LambertW}, where $W(x)$ is Lambert's $W$ function. By the only Proposition in Section 6 of \cite{obaid2013number}, $$2(n+2)^{n-1} = \sum_{i=0}^{n} \binom{n}{i} a(i)a(n-i).$$
So it follows from Lemma \ref{convolutions} that 
\begin{equation} \label{eq:erl}(g(x))^2 = \sum_{n=0}^{\infty} 2(n+2)^{n-1}\frac{x^{n}}{n!}.
\end{equation}
We make the following observations about 
$$ \sum_{i=1}^{n} \binom{n-1}{i-1}(n-i+1)^{(n-i-1)}\cdot i^{i-2}.$$
With the change of variable $j =  i-1,$
$$ \sum_{i=1}^{n} \binom{n-1}{i-1}(n-i+1)^{(n-i-1)} \cdot i^{i-2} = \sum_{j=0}^{n-1} \binom{n-1}{j}(n-j)^{(n-j-2)}(j+1)^{j-1} = 2(n+1)^{n-2},$$
as shown in the proof of the only Proposition in Section 6 of \cite{obaid2013number}. We also observe that
$$ \sum_{i=1}^{n} \binom{n-1}{i-1}(n-i+1)^{(n-i-1)}\cdot i^{i-2} = \sum_{i=1}^{n-3} \binom{n-1}{i-1}(n-i+1)^{(n-i-1)}\cdot i^{i-2} + n^{n-2} + (n-1)^{n-2} + \frac{3}{2}(n-1)(n-2)^{n-3}, $$
so, 
\begin{equation}\label{eq:long} \sum_{i=1}^{n-3} \binom{n-1}{i-1}(n-i+1)^{(n-i-1)}\cdot i^{i-2} =  2(n+1)^{n-2} - n^{n-2} - (n-1)^{n-2} - \frac{3}{2}(n-1)(n-2)^{n-3}.
\end{equation}

As a result we can write the recurrence for $K_{n+1}$ (from Theorem \ref{kn-reccurence}) in the following way ,
\begin{equation} \label{eq:eq2}
 K_{n+1} = 2\cdot2(n+2)^{n-1} -(n+1)^{n-1} - n^{n-1} - \frac{1}{2}(n)(n-1)^{n-2} + \sum_{i=1}^{n-1} \binom{n}{i-1}(n-i)K_{n+1-i} \cdot i^{i-2}.
 \end{equation}
Making the change of variable $j = i-1$ we get, 
$$ K_{n+1} = 2\cdot2(n+2)^{n-1} -(n+1)^{n-1} - n^{n-1} - \frac{1}{2}(n)(n-1)^{n-2} + \sum_{j=0}^{n-2} \binom{n}{j}(n-j-1)K_{n-j} \cdot (j+1)^{j-1}.$$
We will now study the exponential generating function of $K_{n+1}$. To do this we look at the exponential generating function of each of the summands on the right hand side. We have already seen from (\ref{eq:erl}) that
\begin{equation}\label{eq:term1}
 \sum_{n=0}^{\infty} 2(n+2)^{n-1} \frac{x^{n}}{n!} = (g(x))^2
\end{equation}
To deal with the rest of the summands of $K_{n+1}$ in (\ref{eq:eq2}) but the last one, we first re-organise them in the following way using equation (\ref{eq:long}). 
Let $$\phi(n) = \sum_{i=1}^{n-2} \binom{n-1}{i-1}(n-i+1)^{(n-i-1)}\cdot i^{i-2} = 2(n+2)^{n-1}-(n+1)^{n-1}-n^{n-1}-\frac{3}{2}n(n-1)^{n-2}. $$ The change of variable $j=i-1$ gives us
$$\phi(n) = \sum_{j=0}^{n-3} \binom{n}{j}(n-j)^{n-j-2}(j+1)^{j-1}.$$ We have $\phi(n)=0$ for $n=0,1,2$ since the sum is empty for these values of $n$. This further implies that,
$$\sum_{n=0}^{\infty} \phi(n) \frac{x^{n}}{n!} = \sum_{n=2}^{\infty} \phi(n) \frac{x^{n}}{n!}$$
$$ = \sum_{n=2}^{\infty} 2(n+2)^{n-1}\frac{x^{n}}{n!} - \sum_{n=2}^{\infty} (n+1)^{n-1} \frac{x^{n}}{n!} - \sum_{n=2}^{\infty}n^{n-1}\frac{x^{n}}{n!} - \frac{3}{2} \sum_{n=2}^{\infty} n(n-1)^{n-2}\frac{x^{n}}{n!}$$
$$= \left((g(x))^{2} -1 -2x\right) - \sum_{n=2}^{\infty} (n+1)^{n-1} \frac{x^{n}}{n!} - \sum_{n=2}^{\infty}n^{n-1}\frac{x^{n}}{n!} - \frac{3}{2} \sum_{n=2}^{\infty} n(n-1)^{n-2}\frac{x^{n}}{n!},$$ by (\ref{eq:erl}). 
Lemma \ref{LambertW} resolves the second summand. The third summand is resolved by \cite{corless1996lambertw} in Section 2, page 4. Originally this was done in \cite{polya1937kombinatorische}. This has been translated into English; see\cite{polya2012combinatorial}). To resolve the fourth summand we use the fact the exponential generating function is a right index shift and multiplication by $n$ of the 3rd summand. Right index shifting is equivalent to formal integration and by Rule $2^{\prime}$ in Section 2.3 page 41 of \cite{wilf2005generatingfunctionology} multiplication by $n$ is equivalent to differentiating and then multiplying the exponential generating function by $x$ (This is also given on the OEIS A055541). Therefore.
$$\sum_{n=0}^{\infty} \phi(n) \frac{x^{n}}{n!} =  [e^{-2W(-x)}-1-2x] - [e^{-W(-x)}-1-x] - [-W(-x)-x] - \frac{3}{2}[-xW(-x)]$$
$$ = e^{-2W(-x)} -1-2x - e^{-W(-x)} + 1+ x+W(-x)+x+\frac{3}{2}xW(-x)$$
\begin{equation}\label{eq:term2}
= e^{-2W(-x)} - e^{-W(-x)} + W(-x) + \frac{3}{2}xW(-x).
\end{equation}
Now let us study the final summand of (\ref{eq:eq2})  $$\sum_{j=0}^{n-2} \binom{n}{j}(n-j-1)K_{n-j} \cdot (j+1)^{j-1}.$$ Notice that the term $\binom{n}{n-1}(n-(n-1)-1)K_{1}n^{n-2}=0$ and $\binom{n}{n}(n-n-1)K_{0}{n+1}^{n-1}=0$ since $K_{0}=0$. Therefore,
$$\sum_{j=0}^{n-2} \binom{n}{j}(n-j-1)K_{n-j} \cdot (j+1)^{j-1} = \sum_{j=0}^{n} \binom{n}{j}(n-j-1)K_{n-j} \cdot (j+1)^{j-1}.$$
By Lemma \ref{convolutions}, 
$$\sum_{n=0}^{\infty}\left(\sum_{j=0}^{n} \binom{n}{j}(n-j-1)K_{n-j} \cdot (j+1)^{j-1}\right)\frac{x^{n}}{n!} =\left(\sum_{n=0}^{\infty}(n-1)K_{n}\frac{x^{n}}{n!}\right)\left(\sum_{n=0}^{\infty}(n+1)^{n-1}\frac{x^{n}}{n!}\right).$$
By Rule $2^{\prime}$ in Section 2.3 page 41 of \cite{wilf2005generatingfunctionology}
$$\left(\sum_{n=0}^{\infty}(n-1)K_{n}\frac{x^{n}}{n!}\right) = x\frac{d}{dx}\left(\sum_{n=0}^{\infty}K_{n}\frac{x^{n}}{n!}\right)-\left(\sum_{n=0}^{\infty}K_{n}\frac{x^{n}}{n!}\right)=xh^{\prime}(x)-h(x).$$
By Lemma \ref{LambertW}, 
$$\left(\sum_{n=0}^{\infty}(n+1)^{n-1}\frac{x^{n}}{n!}\right) = e^{-W(-x)}.$$
Therefore 
\begin{equation}\label{eq:term3}
\sum_{n=0}^{\infty}\left(\sum_{j=0}^{n} \binom{n}{j}(n-j-1)K_{n-j} \cdot (j+1)^{j-1}\right)\frac{x^{n}}{n!} = (xh^{\prime}(x)-h(x))e^{-W(-x)}.
\end{equation}
By Rule $1^{\prime}$ in Section 2.3 page 41 of \cite{wilf2005generatingfunctionology},
$$\sum_{n=0}^{\infty}K_{n+1}\frac{x^{n}}{n!} = h^{\prime}(x).$$
We now write the exponential generating function of $K_{n+1}$, using the expression of $K_{n+1}$ in (\ref{eq:eq2}) and the exponential generating functions of the summands of $K_{n+1}$ obtained as in (\ref{eq:term1}), (\ref{eq:term2}) and (\ref{eq:term3}). 
$$\sum_{n=0}^{\infty}K_{n+1}\frac{x^{n}}{n!} = 2e^{-2W(-x)} - e^{-W(-x)} + W(-x) + \frac{1}{2}xW(-x)+ (xh^{\prime}(x)-h(x))e^{-W(-x)},$$
so
$$h^{\prime}(x) = 2e^{-2W(-x)} - e^{-W(-x)} + W(-x) + \frac{1}{2}xW(-x)+ (xh^{\prime}(x)-h(x))e^{-W(-x)},$$
Therefore we have the following first order linear ODE,
$$h^{\prime}(x)+ \frac{h(x)e^{-W(-x)}}{(1-xe^{-W(-x)})} = \frac{2e^{-2W(-x)} - e^{-W(-x)} + W(-x) + \frac{1}{2}xW(-x)}{(1-xe^{-W(-x)})}.$$
This ODE is of the form, 
$$h'(x)+Q(x)h(x) = F(x),$$
so we may apply the integrating factor method and give a general solution for $h(x)$,
$$h(x) = e^{-V(x)}\int{V(x)F(x) dx} + C,$$
where $V(X)$ is the integrating factor,
$$V(X) = \int{Q(x) dx} = \int{\frac{e^{-W(-x)}}{1-xe^{-W(-x)}} dx}.$$
Unfortunately, we are unable to evaluate $V(X)$ so we leave $h(x)$ as it is.

\end{proof}
\end{theorem}

\section{Justification} 
In this section we would like to justify why we only look at the four cases above. Our approach to counting the number of complete $\t$-exceptional sequences in the above module categories relied upon Theorems \ref{main1} and \ref{main2}. We also took advantage of the fact that the $\t$-perpendicular categories of indecomposable modules $M$ were of the form $J(M) \cong \mathcal{C} \oplus \mathcal{D}$ with $\mathcal{C} \text{ and } \mathcal{D}$ being module categories in the the two families $\Gamma_{n}^{t}$ or $\Lambda_{n}^{t}$. It is our claim that these four cases, $\Gamma_{n}^{2}, \Gamma_{n}^{n-1}, \Lambda_{n}^{2}, \Lambda_{n}^{n}$ are the only ones were all the $\t$-perpendicular categories $J(M)$ are of this form. In other words, our approach only works on these four cases. 

\begin{proposition} Fix a positive integers $t \geq 3$. For $n \geq t+1$, let $A=\Lambda_{n}^{t}$. Then there exists an $A$-module $M$ such that the $\t$-perpendicular category $J(M)$ is not a direct sum of module categories over algebras of the form $\Lambda_{n^{\prime}}^{t^{\prime}}$ or $\Gamma_{n^{\prime}}^{t^{\prime}}$ for $2 \leq t^{\prime} \leq n^{\prime} < n$.
\begin{proof}
We prove this by counter-example. Set $M = S_{1}$, the simple module at vertex 1 of the quiver $C_{n}$ of $A$. Note that other simple modules also work, but for simplicity we choose $S_1$. The Auslander-Reiten translate of $S_1$ is $\t S_{1} = S_{2}$. Using Proposition \ref{radAssem} and \ref{homAdachi}, we can say that Hom$(X, S_{2}) \neq 0$ if and only if $X=P_{2}/\text{rad}^{l(X)}(P_{2})$ where $l(X)$ is the length of $X$.  It also follows that $P_{2}$ is the only projective with non-zero maps to $S_2$. Therefore all other indecomposable projective modules $P_j$ with $j \neq 2$ are in $^{\perp}(\t S_{1})$, hence they are Ext-projectives in $^{\perp}(\t S_{1})$. By Proposition \ref{ExtProj}, the module $S_1$ is Ext-project in $^{\perp}(\t S_{1})$. We can thus conclude that, $$\mathcal{P}(^{\perp}(\t S_1)) = \bigoplus_{j \neq 2} P_{j} \oplus S_1.$$
By definition the Bongartz completion of $M$ in $\module A$ is $T_{M} = \mathcal{P}(\t S_{1}).$ Let $Q_{n}$ be the following quiver,
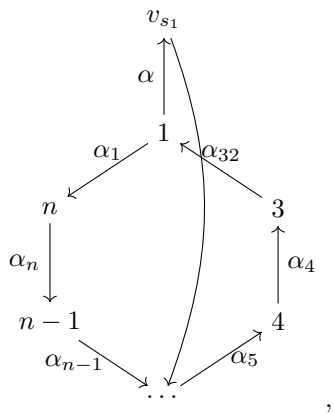
\begin{center}
\begin{tikzpicture}
\node (s1) at (0,5) {$v_{s_{1}}$};
\node (dots) at (0,0) {$\dots$};
\node (4) at (1.5,1) {$4$};
\node (3) at (1.5,2.5) {$3$};
\node (1) at (0,3.5) {$1$};
\node (n) at (-1.5,2.5) {$n$};
\node (n-1) at (-1.5,1) {$n-1$};

\draw[->] (s1) to [out=-70, in =70]  (dots);
\draw[->] (1) -- node[midway, left]{$\alpha$} (s1);
\draw[->] (n) -- node[midway, left]{$\alpha_{n}$} (n-1);
\draw[->] (n-1) --  node[midway, left]{$\alpha_{n-1}$} (dots);
\draw[->] (1) --node[midway,above]{$\alpha_{1}$}  (n);
\draw[->] (dots) --node[midway,right]{$\alpha_{5}$}  (4);
\draw[->] (4) --node[midway, right]{$\alpha_{4}$} (3);
\draw[->] (3) --node[midway, above]{$\alpha_{32}$}  (1);
\end{tikzpicture},
\end{center}
where the vertices labelled $j$ correspond to the projective $P_{j}$ and the vertex $v_{s_{1}}$ corresponds to the simple $S_1$ and the arrows correspond to the irreducible maps between their respective modules. The $\mathbb{F}$-algebra $E_{M} = \text{End}_{A}(T_{M})$ is the path algebra of the quiver modulo relations.  Let $Q_{n}^{v_{s_{1}}}$ be the quiver obtained from $Q_{n}$ by removing the vertex $v_{s_{1}}$ and any arrows incident to $v_{s_{1}}$,
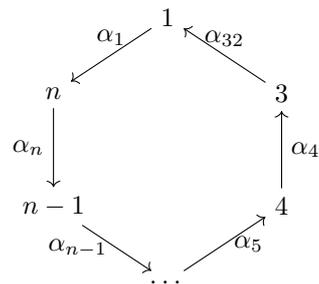
\begin{center}
\begin{tikzpicture}

\node (dots) at (0,0) {$\dots$};
\node (4) at (1.5,1) {$4$};
\node (3) at (1.5,2.5) {$3$};
\node (1) at (0,3.5) {$1$};
\node (n) at (-1.5,2.5) {$n$};
\node (n-1) at (-1.5,1) {$n-1$};

\draw[->] (n) -- node[midway, left]{$\alpha_{n}$} (n-1);
\draw[->] (n-1) --  node[midway, left]{$\alpha_{n-1}$} (dots);
\draw[->] (1) --node[midway,above]{$\alpha_{1}$}  (n);
\draw[->] (dots) --node[midway,right]{$\alpha_{5}$}  (4);
\draw[->] (4) --node[midway, right]{$\alpha_{4}$} (3);
\draw[->] (3) --node[midway, above]{$\alpha_{32}$}  (1);
\end{tikzpicture}
\end{center}
by Lemma \ref{mainlemma}, $D_M = E_{M}/ \langle e_M \rangle$ is the path algebra of the quiver $Q_{n}^{v_{s_{1}}}$ modulo relations.
We have the relation $\alpha_{t+1}\alpha_{t} \dots \alpha_{4}\alpha_{32} = 0$ involving $t-1$ arrows because it corresponds to Hom$_{A}(P_{t+1},P_{2})=0$ since in $ \module \Lambda_{n}^{t}$ the composition of $t$ maps between projectives is 0. However, at the same time we have that the composition of the $t^{\prime}-1$ arrows $\alpha_{1}\alpha_{n} \dots \alpha_{n-(t^{\prime}-3)} \neq 0$ for $2 \leq t^{\prime} \leq t$. Therefore as a module category $J(M)$ cannot be a direct sum of module categories of the form $\module \Lambda_{n^{\prime}}^{t^{\prime}}$  or $\module \Gamma_{n^{\prime}}^{t^{\prime}}$ as required. 
\end{proof}
\end{proposition}

\begin{proposition} Fix a positive integers $t \geq 3$. For $n \geq t+2$, let $A=\Gamma_{n}^{t}$. Then there exists an $A$-module $M$ such that the $\t$-perpendicular category $J(M)$ is not a direct sum of module categories over algebras of the form $\Gamma_{n^{\prime}}^{t^{\prime}}$ or  $\Lambda_{n^{\prime}}^{t^{\prime}}$ for $2 \leq t^{\prime} \leq n^{\prime} < n.$
\begin{proof} 
The argument is similar to that for the previous proposition. We prove this by counter-example. Set $M = S_{1}$, the simple module at vertex 1 of the quiver $A_{n}$ of $A$. The Auslander-Reiten translate of $S_1$ is $\t S_{1} = S_{2}$. By Proposition \ref{radAssem} and \ref{homAdachi}, Hom$(X, S_{2}) \neq 0$ if and only if $X=P_{2}/\text{rad}^{l(X)}(P_{2})$ where $l(X)$ is the length of $X$.  It also follows that $P_{2}$ is the only projective with non-zero maps to $S_2$. Therefore all other indecomposable projective modules $P_j$ with $j \neq 2$ are in $^{\perp}(\t S_{1})$, hence they are Ext-projectives in $^{\perp}(\t S_{1})$. By Proposition \ref{ExtProj}, the module $S_1$ is Ext-project in $^{\perp}(\t S_{1})$ We can thus conclude that, $$\mathcal{P}(\t S_1) = \bigoplus_{j \neq 2} P_{j} \oplus S_1.$$
By definition the Bongartz completion of $M$ in $\module A$ is $T_{M} = \mathcal{P}(\t S_{1}).$ Let $Q_{n}$ be the following quiver,

\begin{center}
\begin{tikzcd}
    v_{s_{1}} & \arrow{l}[swap]{\alpha} 1 & 3 \arrow{l}[swap]{\alpha_{32}} & 4 \arrow{l}[swap]{\alpha_{4}} & \dots \arrow{l}[swap]{\alpha_{5}} & n-1 \arrow{l}[swap]{\alpha_{n-1}} & n  \arrow{l}[swap]{\alpha_{n}}
\end{tikzcd},
\end{center}
where the vertices labelled $j$ correspond to the projective $P_{j}$ and the vertex $v_{s_{1}}$ corresponds to the simple $S_1$ and the arrows correspond to the irreducible maps between their respective modules. The $\mathbb{F}$-algebra $E_{M} = \text{End}_{A}(T_{M})$ is the path algebra of the quiver $Q_{n}$ modulo relations. Let $Q_{n}^{v_{s_{1}}}$ be the quiver obtained from $Q_{n}$ by removing the vertex $v_{s_{1}}$ and any arrows incident to $v_{s_{1}}$,
\begin{center}
\begin{tikzcd}
  1 & 3 \arrow{l}[swap]{\alpha_{32}} & 4 \arrow{l}[swap]{\alpha_{4}} & \dots \arrow{l}[swap]{\alpha_{5}} & n-1 \arrow{l}[swap]{\alpha_{n-1}} & n  \arrow{l}[swap]{\alpha_{n}}
\end{tikzcd}.
\end{center}
by Lemma \ref{mainlemma}, $D_M = E_{M}/ \langle e_M \rangle$ is the path algebra of the quiver $Q_{n}^{v_{s_{1}}}$ modulo relations.
We have the relation $\alpha_{t+1}\alpha_{t} \dots \alpha_{4}\alpha_{32} = 0$ involving $t-1$ arrows  because it corresponds to Hom$_{A}(P_{t+1},P_{2})=0$ since in $\module \Gamma_{n}^{t}$ the composition of $t$ maps between projectives is 0. However, at the same time we have that the composition of the $t^{\prime}-1$ arrows $\alpha_{n}\alpha_{n-1} \dots \alpha_{n-(t^{\prime}-2)} \neq 0$. Therefore as a module category $J(M)$ cannot be a direct sum of module categories of the form $\module \Lambda_{n^{\prime}}^{t^{\prime}}$  or $\module \Gamma_{n^{\prime}}^{t^{\prime}}$ as required. 
\end{proof}
\end{proposition}
So we have shown that our strategy for deriving recurrences for the number of complete $\t$-exceptional sequences over Nakayama algebras only works in the four cases we've studied. However, the statements of Theorems \ref{main1} and \ref{main2} are general enough that a similar strategy may be applied to other algebras, and may prove as effective for counting the $\t$-exceptional sequences for the module categories of those algebras.

\nocite{OEIS}
\nocite{prudnikov1998integrals}
\bibliographystyle{alpha}
	\bibliography{ExceptionalSeqs}
\end{document}